\documentclass[11pt,twoside]{article}
\usepackage{amsmath, amssymb, amsthm}
\usepackage{graphicx} 

\theoremstyle{plain}
\newtheorem*{thmA}{Theorem A}
\newtheorem*{thmB}{Theorem B}

\newtheorem{thm}{Theorem}[section]
\newtheorem{prop}[thm]{Proposition}
\newtheorem{lemma}[thm]{Lemma}
\newtheorem{corollary}[thm]{Corollary}

\theoremstyle{definition}
\newtheorem{defn}[thm]{Definition}

\newtheorem{question}{Question}

\begin{document}

\def\Cal#1{{\cal#1}}
\def\<{\langle}\def\>{\rangle}
\def\what{\widehat}\def\wtil{\widetilde}
\def\Z{{\mathbb Z}}\def\N{{\mathbb N}} \def\C{{\mathbb C}}
\def\Q{{\mathbb Q}}\def\R{{\mathbb R}} \def\E{{\mathbb E}}

\def\Proof{\paragraph{Proof.}}
\def\Remark{\paragraph{Remark.}}
\def\endproof{\hfill$\square$\break\medskip}
\def\noproof{\hfill$\square$\break}
\def\noi{\noindent}

\def\notation{\paragraph{Notation.}}
\def\ackn{\paragraph{Acknowledgement.}}

\let\demph\textbf

\def\diam{\text{diam}}
\def\Lk{\text{Lk}}
\def\ov{\overline}

\def\al{\alpha}                 \def\be{\beta}
\def\ga{\gamma}			\def\Ga{\Gamma}
\def\sig{\sigma}
\def\ep{\epsilon}               \def\varep{\varepsilon}
%

\title{\bf{Two-dimensional Artin groups with CAT(0)\\
 dimension three
\footnote{ The first author acknowledges support from 
the NSF and EPSRC. 
The second author acknowledges the support of a UK EPSRC Research Assistantship, 
and a grant from the Conseil R\'egional de Bourgogne.
Both authors wish to thank the Faculty of Mathematical Studies, 
University of Southampton, for support, conversations and 
hospitality during the preparation of this work.}
}}
 
\author{
\textsc{Noel Brady}\\
\\
\emph{Department of Mathematics, University of Oklahoma,}\\
\emph{Norman, OK 73019, USA}\\
\emph{e-mail:} \texttt{nbrady@math.ou.edu}\\
\\
\textsc{John Crisp}\\
\\
\emph{Laboratoire de Topologie, Universit\'e de Bourgogne,
UMR 5584 du CNRS,}\\
\emph{B.P. 47 870, 21078 Dijon, France}\\
\emph{e-mail:} \texttt{crisp@topolog.u-bourgogne.fr}
}

\date{\emph{(December 7, 2000)}}

\maketitle

\begin{abstract} 
We exhibit 3-generator Artin groups which have finite
$2$-dimensional Eilenberg-Mac~Lane spaces, but which do not act 
properly discontinuously by semi-simple isometries on a 2-dimensional CAT(0) 
complex.  We prove that infinitely many of these 
groups are the fundamental groups of compact, non-positively curved 
3-complexes. These examples show that the 
geometric dimension of a CAT(0) group may be strictly 
less than its CAT(0) dimension.  
\end{abstract}

One way of associating a dimension to a discrete group $G$ is to 
take the minimal dimension of a contractible free $G$-CW complex or, 
equivalently, to take the minimal dimension of a $K(G,1)$ complex. 
This is called the geometric dimension of $G$.  For example, 
groups with torsion are infinite dimensional, and the dimension of 
${\Z}^n$ is equal to $n$. 
Contractible free $G$-CW complexes arise naturally in 
geometric group theory, where the contractibility of the $G$-CW complex 
is often a consequence of non-positive curvature. Specific instances 
include the case 
when $G$ acts  freely and properly discontinuously by isometries 
on a CAT(0) complex, and the case  when a torsion free word-hyperbolic $G$
acts on its Rips complex. 

Restricting attention now to the class of  CAT(0) metric spaces which have the underlying
topology of a CW-complex, we define the {\it CAT(0) dimension} of a group $G$ 
to be the minimal dimension of such a space on which $G$ acts properly discontinuously by 
semi-simple isometries, or $\infty$ if no such action exists.
 (See section~\ref{Sect1} for the definitions of ``CAT(0)'' and ``semi-simple''). 
Clearly the geometric dimension of a torsion free group is 
bounded above by its CAT(0) dimension. 
Examples of groups with finite geometric dimension but infinite CAT(0) dimension are quite easy to find.
Indeed, comparison of translation lengths shows that any polycyclic group which is not virtually abelian fails to act semi-simply on any CAT(0)
space.\footnote{ We thank the referee for pointing out this class of
examples, as well as for several other helpful comments.} 
These include the solvable Baumslag-Solitar groups
$\< a,t\mid tat^{-1}=a^n\>$ for $|n|\neq 1$, 
which have geometric dimension 2. However, it seems that the 
question of \emph{finite} gaps between
the geometric and CAT(0) dimensions is rather more subtle.


In this paper we exhibit an infinite family of Artin groups which all have  
geometric dimension 2, but  have CAT(0) dimension 3. These 
examples illustrate an important point for researchers who endeavour 
to prove that certain families of groups are CAT(0); namely, one should not just 
try to put a non-positively curved metric on a minimal dimensional $K(G,1)$, but should be prepared
to work with higher dimensional $K(G,1)$'s for a given group $G$. 

Our examples come from the following family.
Let $m$, $n$ and $p$ be integers which are all at least 2. We denote by 
$A(m,n,p)$ the three generator Artin group with the following 
presentation
$$
\langle\, a,b,c \; | \; (a,b)_m =  (b,a)_m, \; (b,c)_n =  (c,b)_n, \; 
(a,c)_p =  (c,a)_p \, \rangle
$$
where $(a,b)_m$ denotes the alternating product of length $m$ of $a$'s and 
$b$'s, beginning with $a$. 

There are two key facts which we shall establish about these Artin groups. 
The first is that 3 is a lower bound for the CAT(0) dimension 
in the case where $p=2$ and $m$ and $n$ are both odd.

\begin{thmA} Let $m,n$ be odd integers $\geq 3$. The Artin group $A(m,n,2)$  does not act properly discontinuously by semi-simple isometries on 
any $2$-dimensional CAT(0) complex. 
\end{thmA}

Note that T.~Brady and J.~McCammond \cite{BMcC} have constructed locally 
CAT(0) 2-dimensional Eilenberg-Mac Lane complexes for all $A(m,n,p)$ where 
all $m,n,p\geq 3$, as well as those 3-generator Artin groups with 
fewer than 3 Artin relations. 
In an earlier version of this paper we asked whether any of the groups $A(m,2n,2)$ have CAT(0)
dimension 2. This question has been recently answered P. Hanham \cite{Han} 
who has exhibited compact non-positively curved 2-dimensional piecewise Euclidean
$K(G,1)$ complexes for all of these groups with the exception of the 3-dimensional
groups $A(3,4,2)$, $A(4,4,2)$ and $A(m,2,2)$, $m\geq 2$.

The second key fact is that there are 
explicit 3-dimensional CAT(0) structures for most of the $A(m,n,2)$ groups.

\begin{thmB}
All but finitely many of the 3-generator Artin groups $A(m,n,2)$ are the 
fundamental groups of compact nonpositively curved (i.e: locally CAT(0)) 
3-dimensional piecewise Euclidean complexes. 
\end{thmB}

The methods used to prove Theorem~B can also be used to produce 3-dimensional CAT(0) 
structures for other classes of Artin groups -- see Theorem~\ref{manygens}.

From \cite{CD} we know that the Artin groups $A(m,n,p)$ have geometric dimension 2 precisely when 
$\frac{1}{m} + \frac{1}{n} + \frac{1}{p} \leq 1$. Combining this fact with Theorems~A and~B we get 
that infinitely many of the $A(m,n,2)$ have geometric dimension 2, but CAT(0) dimension 3.
Other examples which exhibit a finite gap between the geometric and CAT(0)
dimensions are given by M. Bridson in \cite{Bri2}. 

We note that our proof of Theorem~A makes no use of the hypothesis that the CAT(0) space in question
has the topology of a CW-complex. Allowing actions on arbitrary CAT(0) spaces and using the covering
dimension of a topological space  leads to the notion of CAT(0) 
dimension which is used in \cite{Bri2} (see \cite{HW} for the theory of covering dimension).
Our arguments are also valid in this context and result in the statement of
Theorem \ref{main4} which is slightly stronger than  Theorem~A.  
There is also a notion of dimension for non-positively curved spaces introduced by B. Kleiner,
in \cite{Kl}, which agrees with the covering dimension on seperable spaces and which is much
better adapted to the study of CAT(0) spaces. We simply remark that Theorem \ref{main4} holds
equally well with respect to Kleiner's dimension. 

Another interesting notion of CAT(0) dimension is obtain by dropping the semi-simplicity assumption and considering all properly discontinuous actions on CAT(0) spaces. In this case examples of a gap between 
CAT(0) geometric dimensions are much harder to find; examples with infinite gap are found 
among the ``random" groups  in Gromov \cite{Gr}--pge.~158. The techniques used in the present paper rely heavily on the assumption of semi-simplicity and we do not know whether Theorem~A is valid without it.   



This paper is organized as follows. In Section~\ref{Sect1} we recall basic 
facts about minsets of isometries of CAT(0) spaces. The dihedral type Artin 
groups are introduced in Section~\ref{Sect2}. Proposition 
\ref{B(a,b)} of Section~\ref{Sect2} is used in the proof of 
Theorem~A to which Section~\ref{Sect4} is devoted. In Section~\ref{Sect3} we 
introduce compact, non-positively curved, 3-complexes for the dihedral type  
Artin groups and use these as building blocks  to construct 
compact, non-positively curved 3-dimensional $K(\pi,1)$ complexes for a 
large family of Artin groups, including those of Theorem~B. 
Section~\ref{Sect5} closes with some remarks and questions that arose at 
various stages during the preparation of this paper.

\section{Minsets and periodic flats in CAT(0) spaces}\label{Sect1}

In this section we recall some standard facts about discrete 
groups of semi-simple isometries of CAT(0) spaces. We shall state the 
flat torus theorem, as well as some facts about the 
structure of minsets of infinite order 
semi-simple isometries. These facts will be used in the proof of the 
Proposition~\ref{B(a,b)}, which is a key ingredient  in the proof 
of Theorem~A. We begin by introducing the notion of a CAT(0) space.

\medskip

A \emph{geodesic}  between points $x$ and $y$ in a metric space $X$ 
is an isometric map of the closed interval $[0,d(x,y)]$ into $X$
such that the endpoints are mapped to $x$ and $y$. We often 
denote the image of such a geodesic by $[x,y]$. 
A \emph{geodesic triangle} $\triangle(x,y,z)$ in a metric space $X$ 
is simply a union of three geodesics $[x,y] \cup [y,z] \cup [z,x]$. 
A \emph{comparison triangle} for $\triangle(x,y,z)$ is a triangle 
$\ov\triangle(\ov x, \ov y, \ov z)$ in the Euclidean plane 
$M^2_1=\E^2$ so that 
$d_X(x,y) \; = \; d_{\E^2}(\ov x,  \ov y)$, 
$d_X(y,z) \; = \; d_{\E^2}(\ov y,  \ov z)$, and 
$d_X(z,x) \; = \; d_{\E^2}(\ov z,  \ov x)$. 
Suppose $p \in [x,y] \subset \triangle (x,y,z)$. A \emph{comparison point} for 
$p$ is a point $\ov p \in \ov\triangle(\ov x, \ov y, \ov z)$ such that $\ov 
p \in [\ov x, \ov y]$ and $d_X(p,x) =  d_{\E^2}(\ov p,  \ov x)$. A geodesic 
triangle $\triangle(x,y,z)$ in a metric space $X$ is said to {\it satisfy 
the CAT(0) inequality} if for all $p,q \in \triangle(x,y,z)$ and comparison 
points $\ov p, \ov q \in \ov\triangle(\ov x, 
\ov y, \ov z)$ we have $d_X(p,q) \leq d_{\E^2}(\ov 
p,  \ov q)$.  A geodesic metric space $X$ (one in which there is a 
geodesic path between any pair of points) is said to be a \emph{CAT(0) 
space} if all geodesic triangles in $X$ satisfy the CAT(0) inequality. 
(One defines the notion of a \emph{CAT(1) space} similarly by comparing every geodesic triangle of total perimeter less than $2\pi$ with its comparison triangle in  the standard $2$-sphere $M^2_1={\mathbb S}^2$.)

The triangle inequality may be usefully rephrased in terms of the Alexandrov (or upper) angle between nontrivial geodesic segments issuing from a common point (see \cite{BH}, p.9, for a definition): $(X,d)$ is a CAT($\kappa$) space ($\kappa =0$ or $1$) if, for every geodesic triangle $\Delta$ with distinct vertices (and of perimeter less than $2\pi$ in the case $\kappa=1$), the Alexandrov angle between any two sides is no greater than the corresponding angle in the comparison triangle $\Delta'$ in $M^2_\kappa$.

A space is \emph{locally} CAT($\kappa$) if every point lies in an open metric ball on which the induced metric is CAT($\kappa$).  The Cartan-Hadamard Theorem (see \cite{BH}) states that a complete metric space $(X,d)$ is CAT(0) if and only if it is locally CAT(0) and simply connected. 

An important source of examples of CAT($\kappa$) spaces come from metric simplicial complexes. These are complexes which are obtained by glueing together metric simplices of constant Riemannian curvature (in our case, either Euclidean or spherical simplices) by isometries of their faces. If one builds $X$ from
a finite number of isometry types of cells (``finitely many shapes'') then a result of Bridson \cite{BH} ensures that the intrinsic pseudometric on $X$ (which extends the metric on each simplex of $X$) is in fact a complete geodesic metric.
   
We will make use of the following criterion for determining whether or not a given piecewise Euclidean complex (with finitely many shapes) is locally CAT(0) (see 
\cite{GdH}~Ch~10~Theorem~15 or \cite{BH}, Theorem II.5.4). A piecewise Euclidean or spherical complex $X$ is said to satisfy the \emph{link condition} if for every vertex $v\in X$ the piecewise spherical link complex $\Lk(v,X)$ is CAT(1). A piecewise Euclidean complex with finitely many shapes is locally CAT(0) if and only if it satisfies the link condition. A piecewise spherical complex (for instance a link complex) is CAT(1) if and only if it satisfies the link condition and contains no closed geodesics (isometrically embedded circles) of length less than $2\pi$. Note that a $1$-dimensional complex (a metric graph) is CAT(1) if and only if every embedded circle has length at least $2\pi$.

In the case of a spherical link complex $\Lk(v,X)$, we distinguish the intrinsic metric $d$ from the metric $d^{(\pi)}(x,y)=\text{min}\{d(x,y),\pi\}$ which measures the Alexandrov angle between geodesic segments emanating from $v$.

 For any point $v$ in a metric space $X$, the \emph{space of directions} $S_v(X)$ is defined as the the set of equivalence classes of geodesic segments emanating from $x$ where two segments are equivalent if the Alexandrov angle between them is $0$ (see \cite{BH}, II.3.18). This is a metric space, with the metric induced by the Alexandrov angle, and corresponds to the link $\Lk(v,X)$ with the metric $d^{(\pi)}$ in the case that $X$ is a piecewise Euclidean complex (or any polyhedron with cells of constant curvature).
\bigskip

We shall work with minsets of isometries of CAT(0) spaces in Sections~\ref{Sect2} and \ref{Sect4}. 
Definition~\ref{minset} through Theorem~\ref{FTT}
establish all the notation  and background facts that we shall use.  

\begin{defn}\label{minset} Let $g$ be an isometry of a 
metric space $X$. The {\it translation length of $g$}, 
denoted by $\ell(g)$, is defined to be 
$$
\ell(g) \; = \; \inf\{d(x,gx) \, | \, x \in X\} \, ,
$$
and the {\it minset of $g$}, denoted by $Min(g)$, is defined to be the 
possibly empty set 
$$
Min(g) \; = \; \{x \in X \, | \, d(x,gx) = \ell(g)\}\,. 
$$
An isometry of a CAT(0) space is called {\it semi-simple} 
if it has a nonempty minset. If $g$ is semi-simple with nonzero translation length then we say that $g$ is \emph{hyperbolic}. In this case any $g$-invariant geodesic line in $X$ shall be called an \emph{axis of $g$} or \emph{$g$-axis}. By \cite{BH} Proposition II.6.2, the $g$-axes all lie in
$Min(g)$. \end{defn}

It is a theorem that a finite order isometry of a CAT(0) space $X$ has a fixed point (see \cite{BH} Corollary II.2.8, for instance). So, if a group acts properly discontinuously by semi-simple isometries on $X$ then it is clear that group elements have infinite order if and only if they are hyperbolic.

We now summarize some fundamental properties of individual isometries and their minsets which may be found in \cite{BH}, Chapter II.6 (see also \cite{Bri1}, Proposition 2.13).
 
Let $g$ be a semi-simple isometry of a CAT(0) space $X$. 
Then $Min(g)$ is always a closed, convex subspace of $X$. 
If $\ell(g)\neq 0$ then $Min(g)$ is isometric to the metric product of 
$\R$ with a CAT(0) space $Y$, where each $\R$-fibre is an axis of $g$.
The element $g$ acts on $Min(g)$ by translating 
along the ${\R}$ factor and fixing the $Y$ factor pointwise. 
Moreover, each isometry of $X$ which commutes with $g$
leaves $Min(g)$ invariant and its restriction to $Min(g)$ also respects
the product structure.

In the 2-dimensional CAT(0) setting the structure of
minsets is particularly simple.
We record this structure in the next Propostion, which follows immediately from the preceding paragraph and Lemma 3.2 of \cite{Bri2}. 

\begin{prop}\label{2Dmins} 
Let  $X$ be a CAT(0) space of covering dimension $2$, and let $g$ be a 
hyperbolic isometry of $X$. Then
\begin{enumerate}
\item  
$Min(g)$ is the product of ${\R}$ and an $\R$-tree $Y$, where each fibre 
${\R} \times \{y\}$ is an axis of $g$. 
\item 
The space of directions $S_v(Min(g))$ of any point $v\in Min(g)$ is the 
orthogonal join of a 0-sphere $S^0 = \{g^+, g^-\}$ (namely, its link in the 
axis of $g$ through $v$) with  a (possibly empty) discrete set. In 
particular, a tripod graph can only be a subspace of $S_v(Min(g))$ if its 
valence three vertex is either $g^+$ or $g^-$.
\end{enumerate}
\end{prop}

Our proof of Theorem~A proceeds by examining a 
particular configuration of periodic flat planes 
in a candidate CAT(0) space for the Artin group. Theorem~\ref{FTT} 
asserts the existence of such periodic flats. 

\begin{thm}[Flat torus theorem: \cite{BH} and 
\cite{Bri1}, Theorem 2.18]\label{FTT}
Let $A$ be a free abelian group of rank $m$ acting properly discontinuously 
by semi-simple isometries on a CAT(0) space $X$.  Then 
$Min(A) = \cap_{a \in A}Min(a)$ is a convex $A$-invariant subspace of $X$ 
which splits as a metric product $Y \times {\E}^n$, where each $a \in A$ 
acts by the identity on the $Y$ factor and by a translation on the $\E^n$ 
factor. Moreover, the quotient of each $\{y\}\times{\E}^n$ by $A$ is an 
$n$-torus. 
\end{thm}

For example, if $\< a,b\>$ is a rank 2 free abelian 
group then $Min(a) \cap Min(b)$ is a union of isometrically embedded flat 
planes on each of which $\langle a, b \rangle$ acts cocompactly as a group
of translations. In the case that $X$ is a CAT(0) space of covering 
dimension $2$ the space $Y$ consists of a single point. (If $Y$ has two 
distinct points then it contains a geodesic segment $I$, and hence $X$ 
contains a closed subset $\E^2\times I$ of covering dimension $3$, a 
contradiction.) Thus there is only one such flat plane. In this case we shall write
$$
\Pi(a,b)= Min(a)\cap Min(b)
$$
for this plane. We make the following important remark that, if there exists 
ANY isometrically embedded copy of $\E^2$ (resp.~$\E^n$) which is invariant 
under $\< a,b\>$ (resp.~$A$) then, by \cite{BH} Proposition II.6.2, this 
plane must lie in $Min(a)\cap Min(b)$ (resp.~$Min(A)$). Thus, when $X$ is 
a $2$-dimensional CAT(0) space, $\Pi(a,b)$ really is the unique 
$\< a,b\>$-invariant flat plane in $X$.

The next lemma lists the possible configurations of periodic flats 
in the minset of a semi-simple isometry in a 2-dimensional 
CAT(0) space.

\begin{lemma}\label{HXBtypes}Let $G$ be a group which acts 
properly discontinuously by semi-simple isometries on 
a  CAT(0) space $X$ of covering dimension $2$.  
Suppose that $\Pi_1 \not= \Pi_2$ are two periodic planes contained in 
$Min(g)$, which correspond to free abelian 
subgroups $\< h_1,g\>$ and $\<h_2,g\>$ of $G$ respectively. 
Then the convex closure of $\Pi_1 \cup \Pi_2$ is one of the 
following types.
\begin{enumerate}
\item ($H$-type or strip-type) In this case $\Pi_1 \cap \Pi_2 = \emptyset$ and 
the convex closure is a quotient space of the disjoint union of 
$\Pi_1$, $\Pi_2$ and a strip  
(product of a closed interval and ${\R}$) where one endpoint of the 
interval times ${\R}$ is identified with a $g$-axis in $\Pi_1$, and the 
other endpoint of the interval times ${\R}$ is identified with 
a $g$-axis in $\Pi_2$.  
\item (X-type) In this case $\Pi_1 \cap \Pi_2$ is a 
$g$-axis and the convex closure is just the union $\Pi_1 \cup \Pi_2$. 
\item (Band-type) In this case $\Pi_1 \cap \Pi_2$ consists of a closed 
interval times ${\R}$, foliated by $g$-axes, and the convex closure 
is again just the union  $\Pi_1 \cup \Pi_2$. 
\end{enumerate}
\end{lemma}

\Proof 
Since $g$ is a hyperbolic isometry of $X$ (it has infinite order in $G$ and acts semi-simply), part (1) of Proposition~\ref{2Dmins} ensures that 
$Min(g)$ is the product of $\R$ and an  $\R$-tree $Y$. Each of the planes 
$\Pi_i$ is contained in $Min(g)$ and so inherits a product structure as 
$\R$ times an axis, $A_i$,  of the $\R$-tree $Y$. 
Since $\Pi_1 \not= \Pi_2$ we have that $A_1 \not= A_2$. 

There are only three possible cases for such axes.  
In the first case, the axes don't intersect. 
The convex closure of these  axes consists of the union of the 
two axes together with a unique segment in $Y$ realizing the 
distance between these axes.  In the second case, the two axes intersect 
in a point in $Y$, and so their convex closure is equal to their union. 
In the third case, the axes intersect along a common interval, and again 
their convex closure is equal to their union.  Now take the product of each 
of these pictures with ${\R}$ to get the three situations described 
in the statement above. 

A priori, there is a fourth possible configuration for $A_1$ and $A_2$. 
They may share a common end. In this case, we can see that 
$\<h_1, h_2\>$ does not act properly discontinuously on $Y$, and 
hence that $\<h_1, h_2, g\>$ does not act properly discontinuously on 
$X$. This contradicts the proper discontinuity of $G$.  
\endproof

\section{The dihedral type Artin groups}\label{Sect2}

In this section we describe the dihedral type Artin groups. We consider a certain pair of rank 2 abelian subgroups of a dihedral type Artin group, and give a proposition which describes how the unique flat planes invariant under these abelian subgroups are arranged in any $2$-dimensional  CAT(0) space on which the whole group acts properly discontinuously by semi-simple isometries.
This information will be used in our proof of Theorem~A.

We begin with some notation which is used to describe the basic Artin relation. 
Given an integer $m \geq 2$ and two letters $a$ and $b$, let $(a,b)_m$ 
denote the alternating product $aba\ldots$ of length $m$ in $a$ and $b$. 
The basic Artin relation is of the form 
\[
(a,b)_m \; = \; (b,a)_m\, .
\]
Given an integer $m \geq 2$, define the dihedral type Artin group 
$A(m)$ to be 
\[
A(m) \; = \; \langle \, a, b \; | \; (a,b)_m = (b,a)_m \, \rangle\, .
\]
One easily sees that $A(2)= \Z^2$. For $m \geq 3$, the group $A(m)$ is  
best understood by factoring out the center $Z$ as follows. 
Let $\Delta = (a,b)_m$.
There are two cases to consider depending on the parity of $m$.

$\bullet$ When $m$ is odd:  $\ Z \; = \;  \langle \Delta^2 \rangle\ $
and $\ A(m)/Z \; =\;  \langle ab \rangle \ast \langle \Delta \rangle 
\; = \; {\Z}_m \ast {\Z}_2 \, .$
In this case we have $\Delta a = b \Delta$ and $\Delta b = a \Delta$.

$\bullet$ When $m$ is even: $\ Z \; = \;  \langle \Delta \rangle\ $
and $\ A(m)/Z \; = \;  \langle a \rangle \ast \langle ab \rangle 
\; =\; {\Z} \ast {\Z}_{m/2}\, .$

\begin{lemma}\label{treeaction}
Suppose that $a$ and $b$ are hyperbolic isometries of an $\R$-tree $T$, such that their product $ab$ has a fixed point. Then the translation axes $\ga_a$ for $a$ and $\ga_b$ for $b$ contain in their intersection a closed interval $J$ whose length is equal to $\text{min}\{\ell(a),\ell(b)\}$. Moreover, $a$ and $b$ translate in opposite directions, as determined on the interval $J$.
(If $\ell(a)\neq\ell(b)$ then in fact $\ga_a\cap\ga_b=J$.)
\end{lemma}

\Proof
We use the fact that, if $g$ is a hyperbolic automorphism of $T$ and $x$ any point in $T$, then $[x,g(x)]$ intersects $\ga_g$ in a subinterval $[x_g,g(x_g)]$ where $x_g$ denotes the projection of $x$ to $\ga_g$. In particular note that $[x_g,g(x_g)]$ is symmetrically placed: $d(x,x_g)=d(g(x),g(x_g))$. 

Since $\ell(ab)=0$, we may suppose that $x\in T$ is fixed by $ab$. Let $y$ denote the point $a^{-1}(x)=b(x)$. Then $[x,y]$ contains intervals $[x_a,y_a]=[x,y]\cap\ga_a$ of length $\ell(a)$, and $[x_b,y_b]=[x,y]\cap\ga_b$ of length $\ell(b)$. Since they are symmetrically placed, one of these two intervals, call it $J$, is contained in the other, and hence in $\ga_a\cap\ga_b$. Clearly $J$ has length $\text{min}\{\ell(a),\ell(b)\}$. Also, it is clear that the direction of translation of $a$ on $\ga_a$ and that of $b$ on $\ga_b$ induce opposing directions on the interval $J$.
\endproof

Consider the group $A(m)$. Let $z_{a,b}$ denote the generator ($\Delta$ or $\Delta^2$) of the centre $Z$ of $A(m)$. 
The following Proposition stems from the existence of rank $2$ abelian subgroups $\< a,z_{a,b}\>$ and $\< b,z_{a,b}\>$ in $A(m)$.

\begin{figure}[ht]
\begin{center}
\includegraphics[width=2.7in]{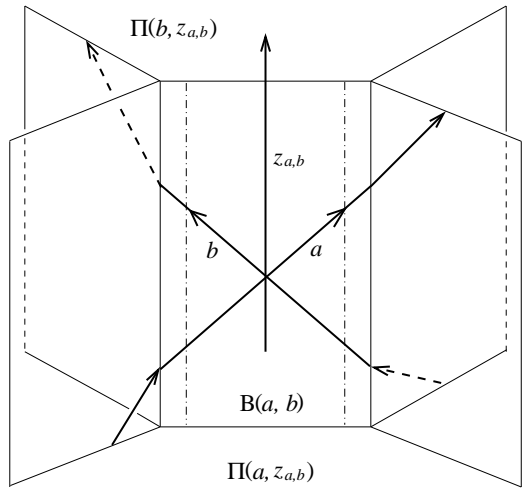}
\end{center}
\caption{Planes $\Pi(a,z_{a,b})$ and $\Pi(b,z_{a,b})$ intersect in a band $B(a,b)$.}\label{Fig1}
\end{figure}

\begin{prop}\label{B(a,b)}
Let $m\geq 3$ and suppose that the dihedral type Artin group
$A(m)$ (with Artin generators $a$ and $b$) 
acts properly discontinuously by semi-simple isometries on a 
CAT(0) space $X$ of covering dimension $2$. 
Let $z_{a,b}$ denote the central element of $A(m)$, and define the set
\[
B(a,b) := \Pi(a,z_{a,b})\cap\Pi(b,z_{a,b}) = Min(a)\cap Min(z_{a,b})\cap Min(b)\,.
\]
Then 
\begin{enumerate}
\item
{$B(a,b)\cong I\times\R$ where $I=[p,q]$ is a nontrivial closed interval and each $\R$-fibre is a $z_{a,b}$-axis.  The $a$-axes intersect the band $B(a,b)$, transversely to the $z_{a,b}$-axis, in a family of parallel line segments which we call \emph{$a$-segments}. Similarly, the $b$-axes intersect $B(a,b)$ in a family of $b$-segments.}
\end{enumerate}

\noindent Moreover, in the case that $m$ is odd, we have the following (See Figure \ref{Fig1}):
\begin{enumerate}
\item [2.] 
The elements $ab, ba \in A(m)$ leave opposite edges of the band  $B(a,b)$ 
invariant. More precisely, $\ell=\{p\}\times\R$ is a $ba$-axis and 
$\ell'=\{q\}\times\R$ is an $ab$-axis. Moreover, $a(\ell)$ is an $ab$-axis 
which lies in $B(a,b)$ and $b(\ell')$ a $ba$-axis also lying in $B(a,b)$.
\item [3.]
For $g,h\in\{ a,b,z_{a,b}\}$ write $\theta(g,h)$ for the angle in $B(a,b)$ between a positively oriented $g$-axis and a positively oriented $h$-axis. Then we have
\[
0<\theta(a,b)=\theta(a,z_{a,b})+\theta(b,z_{a,b})<\pi\ \text{ and } 
\theta(a,z_{a,b})=\theta(b,z_{a,b})\,.
\]
In particular, the $a$-segments are transverse to the $b$-segments.
\end{enumerate}
\end{prop}

\Proof
By Proposition \ref{2Dmins} we have $Min(z_{a,b})\cong\R\times T$ where $T$ 
is an $\R$-tree. Now $A(m)/Z$ acts by isometries on $T$. This action 
is properly discontinuous and 
semi-simple (by  Propositions II.6.10(4) and II.6.9 of \cite{BH}). Since $a$ and $b$ are 
infinite order elements of $A(m)/Z$ they must be hyperbolic, while $ab$ is 
of finite order in $A(m)/Z$ and hence fixes a point in $T$. Thus we may 
apply Lemma \ref{treeaction}. Using the uniqueness of periodic 2-flats in 
the 2-dimensional space $X$, we may identify the plane $\Pi(a,z_{a,b})$ with 
$\gamma_a\times\R$ and $\Pi(b,z_{a,b})$ with $\gamma_b\times\R$. Note that 
$\ga_a\neq\ga_b$ since $A(m)/Z$ is never virtually $\Z$. Thus (1) follows by 
Lemma \ref{treeaction} and Lemma \ref{HXBtypes} (Case 3). In the case that 
$m$ is odd, the element $\Delta$ conjugates $a$ to $b$. It follows that $a$ 
and $b$ have equal translation lengths on both $X$ and $T$. They translate 
in opposite directions in $T$ and it follows from Lemma \ref{treeaction} 
that one endpoint of $I$ is fixed by $ab$, and the other by $ba$. Part (2) 
of the Proposition now follows easily, as does part (3) with the observation 
that that $ba$ translates a positive distance along $\ell$ in the same 
direction as $z_{a,b}$.
\endproof

\section{A lower bound for the CAT(0) dimension of $A(m,n,2)$}\label{Sect4}

In this section we shall prove Theorem~A (Theorem \ref{main4} below) which establishes a lower bound of three for the CAT(0) dimension
of a 3-generator Artin group $A(m,n,2)$ where $m,n \geq 3$ are both odd. It is known \cite{CD} that a 3-generator Artin group $A(m,n,p)$ has geometric dimension 2 for all but the finite number of cases $1/m + 1/n + 1/p>1$ where the geometric dimension is already 3. 
(These are precisely the finite type 3-generator Artin groups). 
These results combined with 
Theorem~\ref{3complexes} yield the result stated in the title of the paper.

\begin{thm}\label{main4}  
Let $G$ be a $3$-generator Artin group of type $(m,n,2)$ where 
$m,n$ are odd integers ($\geq 3$). Then $G$ does not act properly 
discontinuously by semi-simple isometries on a CAT(0) space of covering 
dimension $2$. In particular, the CAT(0) dimension of $G$ is at least $3$, 
and $G$ is not the fundamental group of any compact, non-positively curved 
2-complex. 
\end{thm}

\Remark {(1) As a result of Theorem~\ref{main4}, we see that any Artin group which contains a special subgroup $A(m,n,2)$ with $m,n$ odd, has CAT(0) dimension at least $3$, since the restriction of a semi-simple action to a subgroup always gives a semi-simple action.}

(2) It may be possible to eliminate the hypothesis that $m,n$ be odd.
Its use here is largely confined to the latter part of the proof (Subsection 4.2, Case II).

\Proof We suppose, by way of contradiction, that the Artin group  
\[
A(m,n,2) \; =\; \langle \,a,b,c \; | \; (a,b)_m = (b,a)_m, \, 
(b,c)_n = (c,b)_n, \, [a,c]=1 \, \rangle\, ,\ m,n \text{ odd,} 
\]
acts properly discontinuously by semi-simple isometries on a CAT(0) space 
$X$ of covering dimension $2$. We shall assume throughout that 
$\frac{1}{m}+\frac{1}{n}+\frac{1}{2}\leq 1$ since otherwise $A(m,n,2)$ is a 
finite type Artin group and has geometric dimension $3$. In particular, we 
have $m,n\geq 3$.
 
By Theorem~\ref{FTT} and the fact that $X$ is  2-dimensional, 
each ${\Z} \times {\Z}$ subgroup of $A(m,n,2)$ determines 
a unique isometrically embedded flat plane in $X$. We use the notation $\Pi(x,y)$ for the flat plane $Min(x)\cap Min(y)$ associated to a pair of commuting elements $x$ and $y$. By a well-known result of H. van der Lek \cite{vdL}, the pairs of generators $\{a,b\}$, $\{b,c\}$ and $\{a,c\}$ 
generate the groups $A(m)$, $A(n)$ and $A(2)$ respectively as subgroups of 
$A(m,n,2)$. This gives rise to various $\Z\times\Z$ subgroups of $A(m,n,2)$.
Our proof of Theorem~\ref{main4} involves analyzing the 
configuration of the five flat planes 
\[
\Pi(a,z_{a,b})\, , \quad \Pi(b,z_{a,b})\, , \quad \Pi(b,z_{b,c})\, , \quad \Pi(c,z_{b,c})\, , \quad \Pi(a,c)
\]
corresponding (respectively) to the following rank two free abelian 
subgroups of $A(m,n,2)$ 
\[
\langle\, a,\, z_{a,b}\,\rangle \, , \quad
\langle\, b,\, z_{a,b}\,\rangle \, , \quad
\langle\, b,\, z_{b,c}\,\rangle \, , \quad
\langle\, c,\, z_{b,c}\,\rangle \, , \quad
\langle\, a,\, c\,\rangle\,. 
\]
Here the $z_{s,t}$ denotes the central element in the 2-generator 
Artin subgroup generated by $s$ and $t$. 
Note that $z_{a,c}=ac=ca$ and that
$\Pi(a,z_{a,c}) = \Pi(c,z_{a,c})= \Pi(a,c)$. 
At this point we make no assumption that the five planes just introduced are mutually distinct (even though there may be quite valid group theoretic reasons why they should be).

The first problem we encounter is that we have {\it a priori} no
picture of how all five planes can intersect each other in $X$. 
So we begin by understanding how certain pairs of planes 
arrange themselves in $X$. We then move from this local configuration 
(pairs of planes) to the general configuration of all five planes. 

We begin by considering those pairs of planes 
corresponding to abelian subgroups which share a common central 
element $z_{r,s}$.  For each pair $(r,s)\in\{(a,b), (b,c)\}$, 
Proposition~\ref{B(a,b)} ensures that the periodic flat
planes $\Pi(r,z_{r,s})$ and $\Pi(s,z_{r,s})$, defined above, intersect in 
the nonempty set $B(r,s)=Min(r)\cap Min(z_{r,s})\cap Min(s)$.
Let
\[
\Sigma(r,s)\;  =\; \Pi(r,z_{r,s}) \cup \Pi(s,z_{r,s})
\] 
denote the union of these two flat planes in $X$, and which is in fact a convex subset of $X$ (see Lemma \ref{HXBtypes}).
We also write $\Sigma(a,c) = \Pi(a,c)$, for consistency.  

Next we consider those pairs of planes which correspond to 
abelian subgroups which have an Artin generator in common.  
For the pair of planes $\Pi(a,z_{a,b})$ and $\Pi(a,z_{a,c})$ we define 
$H(a)$ to be the minimal convex subspace of $X$ which contains their union 
as in Lemma~\ref{HXBtypes}. In what follows we will use the terminology introduced in Lemma~\ref{HXBtypes} to distinguish the two cases where these planes are disjoint ($H(a)$ is $H$-type) or intersect ($H(a)$ is $X$-type or $B$-type, or the two planes coincide).
We similarly define $H(b)$ and $H(c)$.

\begin{figure}[ht]
\begin{center}
\includegraphics[width=5in]{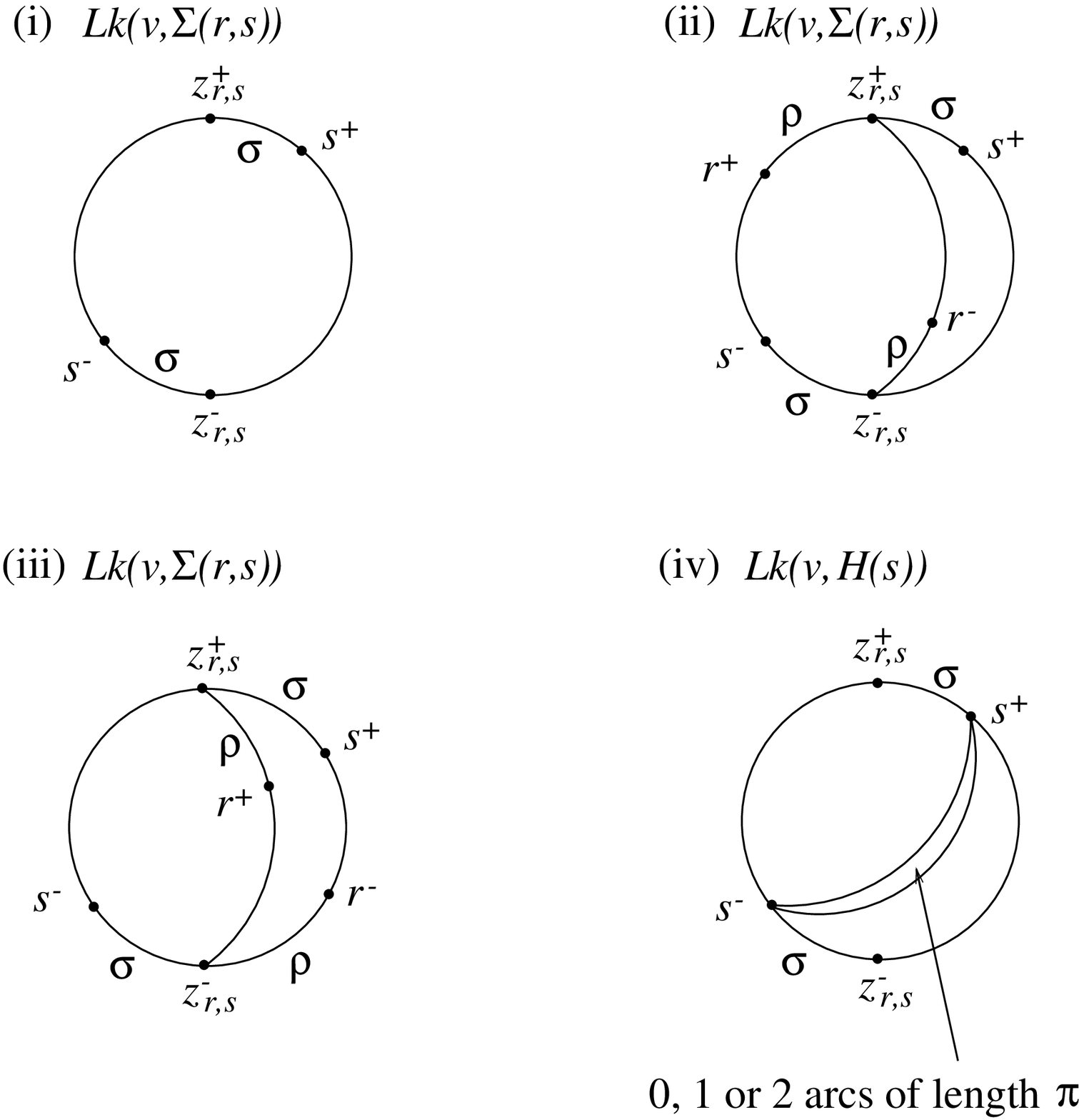}
\end{center}
\caption{Links of a point in $\Sigma(r,s)$ and in $H(s)$.}\label{Fig5}
\end{figure}

We illustrate in Figure~\ref{Fig5} the possible links $\Lk(v,\Sigma(r,s))$ and $\Lk(v,H(s))$ for a point $v \in \Sigma(r,s) \cap H(s)$. 
The angles $\rho=d(r^+,z_{r,s}^+)$
and $\sigma=d(s^+,z_{r,s}^+)$ are indicated in the figure. 
The only constraint on these angles is that $\rho +\sigma < \pi$. This comes from part (3) of Proposition~\ref{B(a,b)}. (At this point we are using the hypothesis that $m,n$ are odd, and in fact $\rho=\sigma$).

The union of the $\Sigma(r,s)$'s and $H(r)$'s forms a (piecewise Euclidean) subspace of $X$ 
which may be viewed as an identification space of the disjoint union 
$\Theta$ of $\Sigma(a,b)$, $\Sigma(b,c)$, $\Pi(a,c)$, $H(a)$, $H(b)$ and 
$H(c)$. However, it is not immediately clear exactly what identifications 
are involved in the map $\Theta\to X$.  We do know, for example, that a 
plane $\Pi(a,z_{a,b})$ in $\Sigma(a,b)$ is identified with a plane in 
$H(a)$, and so on. But there may, a priori, be further identification.  We 
let $Y$ denote the intermediate identification space which results from 
those isometric identifications just mentioned.  Namely, for each of the 
five planes $\Pi(r,z_{r,s})$ one has two isometric embeddings, 
$p:\Pi(r,z_{r,s})\to H(r)\to\Theta$ and 
$q:\Pi(r,z_{r,s})\to\Sigma(r,s)\to\Theta$, and $Y$ is obtained from $\Theta$ 
by identifying the pairs of points $\{ p(x),q(x)\}$ for all 
$x\in\Pi(r,z_{r,s})$ and all $r,s\in\{ a,b,c\}$. 

The space $Y$ may alternatively be described as the union of $H(a),H(b)$ and 
$H(c)$ identified along the subsets $B(a,b)$ (common to $H(a)$ and
$H(b)$), $B(b,c)$ and $B(a,c)=\Pi(a,c)$.  

There is a canonical map $\psi:Y\to X$ induced by the inclusion of each of $H(a)$, $H(b)$ and $H(c)$ into $X$.

The proof of Theorem~\ref{main4} breaks into two cases. In the first case 
$H(a)$ is $H$-type. The case where $H(c)$ is $H$-type is treated similarly. 
We obtain a contradiction in these cases by producing a closed geodesic $\ga$ in $Y$ whose image $\delta=\psi(\ga)$ is either a geodesic 
bigon or a geodesic triangle with angle sum in excess of $\pi$  in the 
CAT(0) space $X$. In the second case neither $H(a)$ nor $H(c)$ are $H$-type 
(so both planes $\Pi(a,z_{a,b})$ and $\Pi(c,z_{b,c})$ intersect $\Pi(a,c)$). 
In this case we  can use the structure of minsets developed in 
Propositions~\ref{2Dmins} and \ref{B(a,b)} to find a fixed point for an 
infinite order element of $A(m,n,2)$. The details of these two cases are 
provided below. This concludes the proof of Theorem~\ref{main4}.  
\endproof

\subsection{\bf Case I : $H(a)$ of $H$-type.} 

\noindent (A similar proof applies in
the case that $H(c)$ is $H$-type).
\medskip

In this case $\Pi(a,c)$ and $\Pi(a,z_{a,b})$ are joined by an infinite strip
$[-1,1]\times\R$ with each line $\{\pm 1\}\times\R$ identified to an
$a$-invariant line in one of the two planes, forming $H(a)$.  Let $\mu$
denote the centre line $\{0\}\times\R$ of the strip and express $H(a)$ as
$T\cup_\mu T'$ where $T=\Pi(a,z_{a,b})\cup([0,1]\times\R)$ and
$T'=\Pi(a,z_{a,c})\cup([-1,0]\times\R)$. Let $\mu,\mu'$ denote the 
preimages of
$\mu$ in $T,T'$ respectively, and let $h:\mu\to\mu'$ denote the
identifying homeomorphism.

Now, cutting $Y$ along $\mu$, we obtain a piecewise Euclidean space

$$
Z= T\cup_{B(a,b)} H(b)\cup_{B(b,c)} H(c)\cup_{B(c,a)} T'\,,
$$

\noindent such that $Y=Z/(x\sim h(x), \hbox{ for all } x\in\mu)$.
\medskip

\begin{lemma}\label{ZYCAT(0)}
The space $Z$, and hence also $Y$, is locally CAT(0).
Moreover, $Z$ is contractible.
\end{lemma}

\Proof We use the link condition. 
Every point in $Z$ has a link which is may be built up, in
a finite number of steps, from a $2\pi$-circle by attaching at each step
an arc of length $\pi$ between antipodal points on a circle of length
$2\pi$ in the link obtained at the previous step. (See Figure~\ref{Fig5}).
The result of such a process will always be a CAT(1) link.
Thus $Z$ is locally CAT(0). 
Links of points in $Y$ are as for points in $Z$ with the exception of
points along $\mu$ where the link is always a $2\pi$-circle.  So $Y$ is
also locally CAT(0).  Finally, one may easily construct a homotopy
retraction of $T$ onto $B(a,b)$, then of $H(b)$ onto $B(b,c)$, then of
$H(c)$ onto $B(c,a)$, then of $T'$ to a point, showing that $Z$ is a
contractible space. 
 
\endproof

Note that $Z$ may be thought of as a fundamental domain for the action
of $\pi_1(Y)\cong\Z$ on the universal cover $\widetilde Y$ of $Y$. 
Moreover, by the Cartan-Hadamard Theorem (\cite{BH}, Theorem II.4.1),
both $\wtil Y$ and $Z$ are CAT(0) spaces.

\begin{lemma}\label{geoloop} There exists a closed 
geodesic $\gamma$ in $Y$ which lifts to a geodesic path
$\wtil\gamma :[0,1]\to Z$  with the following properties.
\begin{description}
\item{(i)} $\wtil\gamma(0)\in\mu$ and 
$\wtil\gamma(1)=h(\wtil\gamma(0))\in\mu'$,
\item{(ii)} there exist $0 \leq t_{a,b}\leq t_{b,c}\leq t_{c,a}\leq 1$ such 
that 
$$
\wtil\gamma(t_{a,b})\in B(a,b)\,,\quad \quad 
\wtil\gamma(t_{b,c})\in B(b,c)\,,\quad \quad 
\wtil\gamma(t_{c,a})\in B(c,a)\,, 
$$
and
$$
\wtil\gamma([0,t_{a,b}])\in T\,,\quad \quad  
\wtil\gamma([t_{a,b},t_{b,c}])\in H(b)\,,\quad \quad 
\wtil\gamma([t_{b,c},t_{c,a}])\in H(c)\,,\quad \quad 
\wtil\gamma([t_{c,a},1])\in T'\,. 
$$
\end{description}
\end{lemma}

\Proof
The space $Y$ may be triangulated, and hence $\wtil Y$ equivariantly  triangulated, into a piecewise Euclidean complex with finitely many
isometry types of simplexes.
Thus, by \cite{BH}, II.6.6 Exercise (2),
every simplicial isometry of $\wtil Y$ is semi-simple. 
In particular, the deck transformation $g:\wtil Y\to\wtil Y$
associated to a generator of 
$\pi_1(Y)\cong\Z$ is a semi-simple with $\ell(g)\neq 0$.
Thus $Min(g)$ contains a $g$-axis whose image in $Y$ is a closed geodesic $\gamma$ with the desired properties. For (i), note that $h=g\mid_\mu$.
Part (ii) follows from the definition of $Z$ and is true for any geodesic from a point on $\mu$ to a point on $\mu'$. We note that
$T$, $H(b)$, $H(c)$ and $T'$ are each closed convex subsets of $Z$.
 
\endproof

We now consider the closed geodesic $\gamma:S^1\to Y$ of Lemma \ref{geoloop}, where $S^1$ is decomposed into a union of closed intervals
$I_a=[t_{c,a},1]\cup [0,t_{a,b}]$, $I_b=[t_{a,b},t_{b,c}]$
and $I_c=[t_{b,c},t_{c,a}]$.
Let $\delta$ denote the closed path in $X$
obtained by composing $\gamma$ with the map $\psi:Y\to X$. 
For each $s\in\{ a,b,c\}$, let $\delta_s$
denote the path obtained by restricting $\delta$ to the interval $I_s$.
Note that the restriction of $\gamma$ to $I_s$ is a locally geodesic path
in $Y$ with image in $H(s)$ and hence defines a geodesic in the CAT(0) space $H(s)$.
Since each of $H(a),H(b)$ and $H(c)$ is a convex subspace of $X$, it now
follows that each of $\delta_a$, $\delta_b$ and $\delta_c$ is a geodesic
path in $X$. 
Thus $\delta$ is a (possibly degenerate) geodesic triangle
in $X$ with sides $\delta_a,\delta_b,\delta_c$ and corners $\delta(t_{a,b}), \delta(t_{b,c}), \delta(t_{c,a})$.

For $r,s\in\{ a,b,c\}$, we say that $t_{r,s}$ is a {\it non-degenerate
corner point} of $\delta$ if both $\delta_r$ and $\delta_s$ are sides
of nonzero length.

Note that at least $\delta_a$ is a side of nonzero length (it must cross the strip $[-1,1]\times\R$ in $H(a)$). Thus we cannot have $t_{a,b}=t_{b,c}=t_{c,a}$ for then $\delta_a$ would be a 
geodesic in $X$ of nonzero length which starts and finishes at the same point, a contradiction. 
 
Therefore, we are in one of the following cases:

\begin{enumerate}
\item 
$\delta$ is a non-degenerate triangle ($t_{a,b}$, $t_{b,c}$ and $t_{c,a}$ are distinct non-degenerate corner points),

\item 
$\delta$ is a bigon with one non-degenerate corner point
whose opposite side has zero length. 

\end{enumerate}

\noindent But it follows from Lemma \ref{corners} below that we now have, in case (1), 
a geodesic triangle with angle sum greater than $\pi$, and in case (2), 
a bigon with at least one nonzero angle.
Both of these outcomes contradict $X$ being a CAT(0) 
space, which completes the reduction of our proof to Case II.

\medskip

Before proving Lemma \ref{corners} we need to establish the following Lemma.

\begin{lemma}\label{anglebound}
Suppose that $f:U\to V$ is a locally isometric map from a geodesic space $U$ to a CAT(1) space $V$. Suppose that $p,q\in U$ are such that $d_U(p,q)\geq \pi$. Then
\[
d_V(f(p),f(q))\geq 2\pi - d_U(p,q)\,.
\]
\end{lemma}

\Proof
Let $[p,q]$ denote a geodesic segment in $U$ from $p$ to $q$, and let $r$ be 
the midpoint of $[p,q]$. We may as well suppose that $d_U(p,q)<2\pi$, in 
which case $[p,q]$ is the union of geodesic segments $[p,r]$ and $[r,q]$ 
each of length $\frac{1}{2}d_U(p,q)$ which is less than $\pi$. 
Since $f$ is locally isometric and $V$ is CAT(1), 
these map to (unique) geodesic 
segments  $[f(p),f(r)]=f([p,r])$ and $[f(r),f(q)]=f([r,q])$ in $V$. 
Moreover, since $f$ is a local isometry at the point $r$, and $[p,q]$ is 
geodesic, we also have that the Alexandrov angle at $f(r)$ between these 
geodesic segments is $\pi$.

Now suppose, by way of contradiction, that $d_V(f(p),f(q))< 2\pi-d_U(p,q)$.
Then, since $d_U(p,q)>\pi$, we have $d_V(f(p),f(q))< \pi$ and therefore
a unique geodesic segment $[f(p),f(q)]$ in $V$. We now have a triangle 
$\Delta=\Delta(f(p),f(q),f(r))$ in $V$ with perimeter $d_V(f(p),f(q))+ 
d_U(p,q)< 2\pi$. But then any comparison triangle in $\mathbb{S}^2$ for 
$\Delta$ has  all angles strictly less than $\pi$. Thus $\Delta$, having an 
angle of $\pi$, fails the CAT(1) inequality (\cite{BH}, Proposition 
II.1.7(4)), a contradiction. 
\endproof

\begin{lemma}\label{corners}
Suppose that, for some $r,s\in\{a,b,c\}$, the point
$t_{r,s}$ is a non-degenerate corner point of $\delta$.  
Then the geodesic paths $\delta_r$ and $\delta_s$ form an angle 
$\theta_{r,s}$ at $\delta(t_{r,s})$ in $X$ where $\theta_{a,c}\geq\frac{\pi}{2}$ and
$\theta_{a,b},\theta_{b,c}\geq\frac{\pi}{3}$. 
\end{lemma}

\Proof Let $v=\delta(t_{r,s})$. Then $v$ lies in each of the closed convex subspaces $H(r)$, $H(s)$ and $\Sigma(r,s)$ of $X$.
We now define the metric graph $L$ to be the identification space 
$$
L = \Lk(v,H(r))\cup \Lk(v,\Sigma(r,s))\cup\Lk(v,H(s))\, / \sim
$$
where $\Lk(v,H(r))$ and $\Lk(v,\Sigma(r,s))$ are identified 
along the common circular link observed for $v$ in $\Pi(r,z_{r,s})$, and
$\Lk(v,H(s))$, $\Lk(v,\Sigma(r,s))$ are identified along the circular
link of $v$ in $\Pi(s,z_{r,s})$. Note that $L$ is a geodesic space, and may 
be viewed as  a sub-link of  the 
link $\Lk(\gamma(t_{r,s}),Y)$ where $\gamma$ is the closed geodesic of 
Lemma~\ref{geoloop}.  

We now make use of the space of directions $S_v(X)$ which may be defined for any point $v$ in a metric space $X$. 
A theorem of Nikolaev \cite{Ni} (see also \cite{BH}) states that if $X$ is a CAT(0) space then the metric completion $S_v(X)'$ of $S_v(X)$ is a CAT(1) space. (Any space is isometrically embedded in its metric completion).

\begin{figure}[ht]
\begin{center}
\includegraphics[width=5.4in]{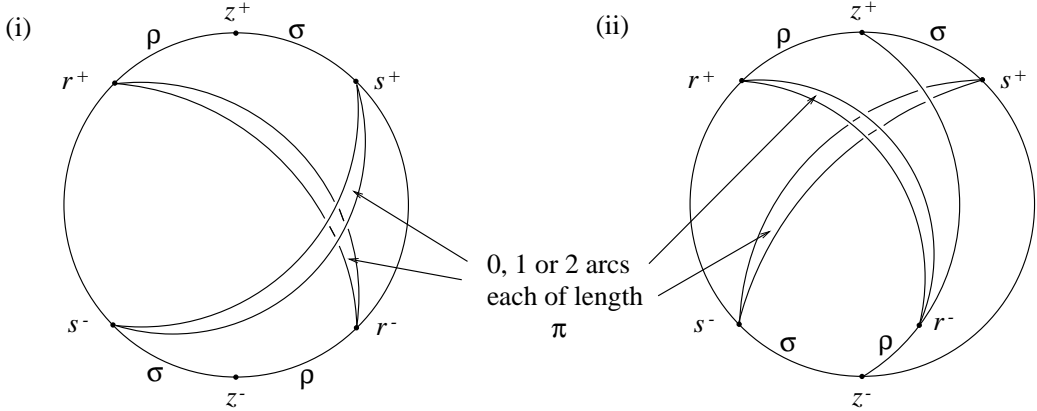}
\end{center}
\caption{Possible forms of the metric graph $L$.}\label{Fig6}
\end{figure}

The inclusions of each subspace into $X$ together induce a map 
$f:L\to S_v(X)\subset S_v(X)'$, which we observe to be locally
isometric as follows. 
A neighbourhood of any point of $L$ lies wholly in the subgraph
$\Lk(v,R)$ where the subspace $R$ is one of $H(r),H(s)$ or $\Sigma(r,s)$.
(For this assertion we are using the fact that $\rho+\sigma < \pi$ so that $r^+$ does not coincide with $s^-$ or $s^+$ with $r^-$ in $L$. See Figure \ref{Fig6} for the possible forms of $L$.)
Since each of these subspaces $R$ embeds in $X$ as a convex subspace,
$f$ restricts to a locally isometric embedding on $\Lk(v,R)$. Hence any point in $L$ possesses a neighbourhood which is embedded isometrically in $S_v(X)'$ by $f$.

We wish to measure the Alexandrov 
angle between $\delta_r$ and $\delta_s$ at $v$. 
Let $x \in \Lk(v,H(r))$ and $y\in\Lk(v,H(s))$ be such that $\delta_r$ 
determines the point   $f(x) \in S_v(X)$ and $\delta_s$ determines the point 
$f(y) \in S_v(X)$. Then the  two sides $\delta_r$ and $\delta_s$  form
an angle at $v$ given by $\theta_{r,s}=d^{(\pi)}_{S_v(X)}(f(x),f(y))$. 
Note that $d_L(x,y)\geq \pi$, since the fact that $\ga$ is a geodesic implies that $d^{(\pi)}(x,y)=\pi$ when viewed in $\Lk(\ga(t_{r,s}),Y)$.
By Lemma \ref{anglebound} we now have $\theta_{r,s}\geq 2\pi-d_L(x,y)\,.$
The result now follows by Lemma \ref{diam(L)} which shows that 
$d_L(x,y)\leq\frac{5\pi}{3}$, and $\leq\frac{3\pi}{2}$ in the case 
$m_{r,s}=2$.
\endproof

\begin{lemma}\label{diam(L)} Let $L$ be the metric graph obtained from 
$\Lk(v, H(r))$, $\Lk(v,\Sigma(r,s))$ and $\Lk(v,H(s))$ as in  
the proof of Lemma~\ref{corners}. 
Define $\diam(L)=sup\{d_L(u,v):u,v\in L\}$.
\begin{description}
\item{(1)} If $m_{r,s}\geq 3$, then 
$\diam(L)\leq \frac{5\pi}{3}$, and 
\item{(2)} if $m_{r,s}=2$, then $\diam(L)\leq\frac{3\pi}{2}$. 
\end{description}
\end{lemma}

\Proof
In view of the possibilities illustrated in Figure~\ref{Fig5}, $L$ takes the
form of one of the links shown in Figure~\ref{Fig6}(i) if $m_{r,s}=2$, that is if $\{ r,s\}=\{a,c\}$,
and otherwise any of Figures~\ref{Fig6}(i) or \ref{Fig6}(ii), or \ref{Fig6}(ii) with the roles of $r$ and $s$ reversed.

Each of the links illustrated in Figure~\ref{Fig6} consists of
$\Lk(v,\Sigma(r,s))$ with zero, one or two arcs of length $\pi$
connecting between $r^+$ and $r^-$, and similarly between $s^+$ and
$s^-$.  We shall call these additional arcs $H(r)$-arcs and $H(s)$-arcs
respectively. 

We consider first the graphs illustrated in Figure~\ref{Fig6}(i). This includes the case where $m_{r,s}=2$. 
In fact if a graph $L$ is formed from a $2\pi$ circle by adding arcs of 
length  $\pi$ onto at most two `great 0-circles' (pairs of antipodal points), 
then the diameter of $L$ is at most $\pi+\phi$ where $\phi$ denotes the 
distance between the two great $0$-circles. One can easily see this by 
producing a family $\{c_i\}$ of closed geodesics of length either $2\pi$ 
or $2(\pi+\phi)$ such that every pair of points in $L$ lie on a 
common $c_i$. In particular $\diam(L)\leq \frac{3\pi }{2}$.

Now consider the links illustrated in Figure~\ref{Fig6}(ii).
Take any $p,q\in L$. Note that by the reasoning just given, the 
union of $\Lk(v,\Sigma(r,s))$ with just the $H(r)$-arcs (or with just the $H(s)$-arcs) 
has diameter at most $\frac{3\pi}{2}$. Thus we may as well assume 
that $p$ lies in an $H(r)$-arc and $q$ lies on an $H(s)$-arc.

Now consider the distances in $\Lk(v,\Sigma(r,s))$ between either of
the points $r^\pm$ and either of $s^\pm$ which may be displayed in a
matrix as follows:
$$
\begin{pmatrix}
d(r^+,s^+) &d(r^+,s^-)\\
d(r^-,s^+) &d(r^-,s^-)
\end{pmatrix}
\; = \; \begin{pmatrix}\rho+\sigma &\pi-(\rho+\sigma)\\
\pi-|\rho-\sigma| &\rho+\sigma
\end{pmatrix}\,.
$$ 

We describe four locally geodesic paths from $p$ to $q$ which we denote
$\tau(\epsilon,\nu)$, for $\epsilon,\nu\in\{+,-\}$ as follows. The path $\tau(\epsilon,\nu)$ passes through
$r^\epsilon$ and $s^\nu$ and follows in part the unique geodesic in
$\Lk(v,\Sigma(r,s))$ of length $d(r^\epsilon,s^\nu)$.  The
paths $\tau(+,+)$ and $\tau(-,-)$ together form a closed geodesic in
$L$ of total length $2\pi +d(r^+,s^+)+d(r^-,s^-)$ which is an upper
bound for $2d(p,q)$.  Similarly with $\tau(+,-)$ and $\tau(-,+)$. 
This gives the following upper bound  for $d(p,q)$ (which is, in fact,
precisely the diameter of $L$):

\[
\begin{aligned}
d(p,q)&\leq \pi + \frac{1}{2}
\hbox{min}\{d(r^+,s^+)+d(r^-,s^-),d(r^+,s^-)+d(r^-,s^+)  \}\\
	&=\pi +
\hbox{min}\{\rho+\sigma,\pi- \frac{1}{2}(\rho+\sigma+|\rho-\sigma|) \}\\
&=\pi +\hbox{min}\{\rho +\sigma,\pi- \hbox{max}\{\rho,\sigma\}\, \}\\
&=\pi +\hbox{min}\{\rho +\sigma, \pi -\rho, \pi -\sigma \} \,.
\end{aligned}
\]

\noindent This upper bound is greatest when
$\rho=\sigma=\frac{\pi}{3}$. Thus we always have
$\diam(L)\leq\frac{5\pi}{3}$.

\endproof

\subsection{\bf Case II : Both $\Pi(a,z_{a,b})$ and
$\Pi(c,z_{b,c})$ intersect $\Pi(a,c)$.} 
\medskip

We obtain a contradiction in this case by proving that the set 
\[
R := Min(a)\cap Min(b)\cap Min(c) = Min(b)\cap \Pi(a,c)
\]

\noindent contains a point which is fixed by an infinite order group element. In order to give a more concrete description we first 
consider the following subsets of $R$. Define
\[
\begin{aligned}
R(a,b) &:= R\cap Min(z_{a,b}) = R\cap\Pi(b,z_{a,b})\ \hbox{ and,}\\
\ &\text{} \\
R(b,c) &:= R\cap Min(z_{b,c}) = R\cap\Pi(b,z_{b,c})\,.
\end{aligned}
\]

\begin{figure}[ht]
\begin{center}
\includegraphics[width=4.2in]{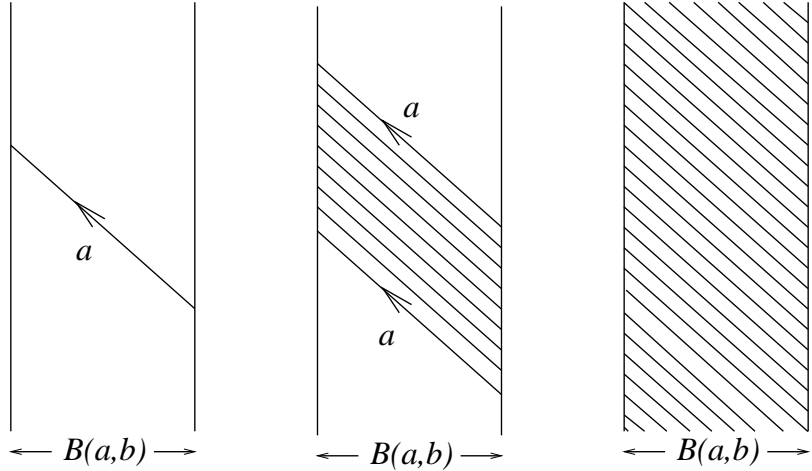}
\end{center}
\caption{$R(a,b)$ as a subset of $B(a,b)$}\label{Fig7}
\end{figure}

\noindent Note that $R(a,b)=B(a,b)\cap\Pi(a,c)$. In fact, if we let $S(a)$ denote the nonempty subset $\Pi(a,z_{a,b})\cap\Pi(a,c)$ of $H(a)$, then $R(a,b)$ should be thought of as the part of $B(a,b)$ which intersects $S(a)$. Now $S(a)$ is a union of $a$-axes: either a single $a$-axis ($H(a)$ is $X$-type), an interval's worth $I\times\R$ ($H(a)$ is $B$-type), or the whole of $\Pi(a,z_{a,b})$. Each $a$-axis (of $\Pi(a,z_{a,b})$) intersects $B(a,b)$ transversely in a closed interval, which we call an \emph{$a$-segment}.
The set $R(a,b)$ is then a closed convex union of (parallel) $a$-segments: either a single $a$-segment, a parallelogram bounded on two sides by 
$a$-segments and on the other sides by the sides of $B(a,b)$, or just the whole of $B(a,b)$ (see Figure \ref{Fig7}).

We have similar possible descriptions for the set $R(b,c)$ as a union of $c$-segments lying in the band $B(b,c)$. We now have the following description of $R$:

\begin{lemma}\label{regionR}
The set $R$ is a 2-dimensional closed convex subset of $\Pi(a,c)$ and $R=R(a, b)= R(b,c)$. 
\end{lemma}

\Proof Clearly, as an intersection of minsets, $R$ is a closed convex 
subset of $\Pi(a,c)$ which contains both $R(a,b)$ and $R(b,c)$.
In particular $R$ contains the convex hull of the union of
an $a$-segment from $R(a,b)$ and a $c$-segment from $R(b,c)$. 
Since these segments lie in different directions in the plane $\Pi(a,c)$, it follows that $R$ is a $2$-dimensional (possibly infinite) closed convex subset of $\Pi(a,c)$.
  
If $R \not= R(a,b)$ then one can find a point $x$ of $R(a,b)$ which 
lies in the closure of $R\setminus R(a,b)$.
Thus $\Lk(x,R)$, which is closed connected subset of ${\mathbb S}^1$,
contains $\Lk(x,R(a,b))$ as a proper closed subset. In particular, 
$\Lk(x,R)$ contains a nontrivial closed arc which intersects $\Lk(x,R(a,b))$ 
at a single point, necessarily one of $a^\pm$ or $z_{a,b}^\pm$.

Note that $R\cup\Pi(b,z_{a,b})$ is a convex subset of $Min(b)$ where $R$ intersects $\Pi(b,z_{a,b})$ along the set $R(a,b)$. Considering the contributions from $R$ and $\Pi(b,z_{a,b})$, we see that 
$S_x(Min(b))$ contains a tripod graph whose 
valence three vertex is one of $a^\pm$ or $z_{a,b}^\pm$. 
But this contradicts  Proposition~\ref{2Dmins}(2) since, by 
Proposition~\ref{B(a,b)} (1) and (3), $b$-axes are transverse to both 
$a$-axes and to $z_{a,b}$-axes. Thus $R = R(a,b)$. 
Similarly, one can show that $R = R(b,c)$.  
\endproof

There are now only two possibilities for the region $R$.
Either $B(a,b)$ and $B(b,c)$ coincide and $R=B(a,b)=B(b,c)$, or $R$ is a nondegenerate parallelogram bounded by a pair of $a$-segments and a pair of $c$-segments. We treat the cases separately. In both cases we apply Proposition \ref{B(a,b)} (2) to each of $B(a,b)$ and $B(b,c)$, thus using the fact that both $m$ and $n$ are odd.
\medskip

\noi{\bf (i) $B(a,b)$ and $B(b,c)$ coincide in $R$.}
\smallskip

\begin{figure}[hb]
\begin{center}
\includegraphics[width=6.1in]{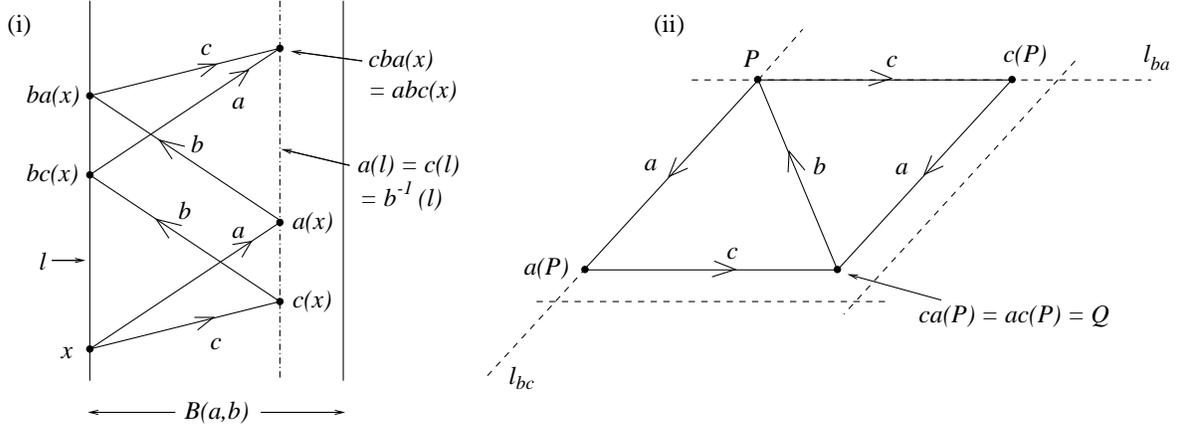}
\end{center}
\caption{Finding a fixed point in $R$}\label{Fig8}
\end{figure}

In this case the $ba$-invariant edge of $B(a,b)$, which we shall denote 
$\ell$, must coincide with the $bc$-invariant edge of $B(b,c)$.
Moreover, since $ba(\ell)=bc(\ell)=\ell$, we have $a(\ell)=c(\ell)=b^{-1}(\ell)$, another line in $R$ parallel to $\ell$. By considering translation lengths in the direction of these lines one now sees that 
$abc(p)=cba(p)$, for all $p\in\ell$. (See Figure~\ref{Fig8}(i)).
Thus every point of $\ell$ is fixed
by the element $g=a^{-1}b^{-1}c^{-1}abc\in A(m,n,2)$. There is a 
canonical map of $A(m,n,2)$ onto the Coxeter group 
\[
W(m,n,2)=\langle a,b,c \mid a^2=b^2=c^2=1, (ab)^m=(bc)^n=(ac)^2=1 \rangle     
\]

\noindent which takes $g$ to the element $(abc)^2$.
When $1/m+1/n+1/2\leq 1$ the group $W(m,n,2)$ is isomorphic to the group generated by reflections in the sides of a Euclidean or hyperbolic triangle with angles $\pi/m$, $\pi/n$ and $\pi/2$, and the element $abc$ is known to be an element of infinite order in this group. Therefore $g$ must be an element of infinite order, giving our contradiction.
\medskip

{\bf (ii) $R$ is a parallelogram.}
\smallskip

We use the fact that the lengths of the $a$-segments and $c$-segments which
bound $R$ in this case are respectively greater than or equal to the translation
lengths $\ell(a)$ and $\ell(c)$. Thus there is a (unique) vertex $P$ of
$R$ for which both $a(P)$ and $c(P)$ lie in $R$. Let $Q=ac(P)=ca(P)$. 
Using the fact that $R=R(a,b)$ and recalling Proposition~\ref{B(a,b)}, we observe that 
the interval segment $[P,c(P)]$ lies on a $ba$-axis $\ell_{ba}$. 
Similarly, by viewing $R$ as $R(b,c)$, the segment $[P,a(P)]$ lies on a 
$bc$-axis $\ell_{bc}$. 

Note that the segment $[P,a(P)]$ (which lies on an $a$-axis) meets 
$\ell_{ba}$ transversely when viewed in the flat plane $\Pi(a,z_{a,b})$. It 
follows that $\ell_{ba}$ and $\ell_{bc}$ intersect transversely at $P$,
and in particular that $\ell_{ba}\cap\ell_{bc}=\{ P\}$.
  
Now, since $c(P)\in\ell_{ba}$ and $\ell_{ba}$ is $ba$-invariant, it follows
that $b(Q)=bac(P)$ also lies in $\ell_{ba}$. But, by the same reasoning, $b(Q)$ lies in $\ell_{bc}$. Therefore we must have $b(Q)=P$. Now it follows that $abc$ fixes the point $a(P)$ (see Figure~\ref{Fig8}(ii)) and, as was observed above, $abc$ is an element of infinite order in $A(m,n,2)$.
Thus we again have a contradiction.
\medskip

This completes the proof of Case II, and of Theorem~\ref{main4}.

\section{Artin groups with $3$-dimensional 
non-positive curvature}\label{Sect3}

In this section we produce 3-dimensional CAT(0) structures for 
a large class of Artin groups, including all 
but finitely many of the 3-generator Artin groups, thus establishing Theorem~B of the Introduction (see Theorem \ref{3complexes} and Corollary \ref{3gens} below).

\subsection{The building blocks}

We begin by constructing a 3-dimensional locally CAT(0) Eilenberg-Mac Lane complex $X_m$  for each dihedral type Artin group 
$$
A(m) \; = \; \langle\,  a,b \; | \; (a,b)_m = (b,a)_m \, \rangle\,, \hbox{for } m\geq 3 .
$$
This complex, at least for the trefoil knot group $A(3)$,  was previously known to a number of researchers 
in the field. It was described to us by B.H. Bowditch who attributed it to Thurston. A description of it is also given by Bestvina \cite{Be}. 

Start with a regular euclidean $m$-gon, $M$, centered on the origin 
in the $xy$-plane. We shall first describe $X_m$ as a certain 
quotient space of the infinite $m$-gon prism, $M\times\R\subset\R^3$. Let $\phi: M \times\R\to M\times\R$ 
be the homeomorphism which translates a positive distance vertically (in $z$-direction) 
and rotates through $(2\pi/m)$ about the $z$-axis. 
Let $E$ denote the union of the vertical edges 
of $M \times\R$; that is, the union of 
edges of the form $\{v\}\times\R$ 
where $v$ is a vertex of $M$. The 3-complex, $X_m$, is defined to be the 
quotient of  $M \times\R$ by the maps $\phi|_{E}$ and  
$\phi^2$.  

\begin{figure}[ht]
\begin{center}
\includegraphics[width=1.5in]{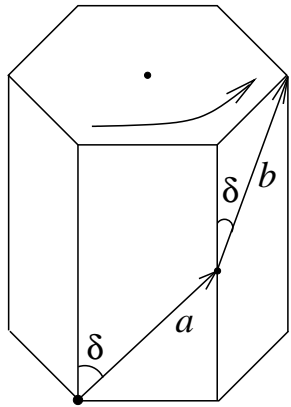}
\end{center}
\caption{The 3-complex for the dihedral type Artin group $A(m)$.}\label{Fig2}
\end{figure}

The universal cover of $X_m$ consists of an infinite collection of 
copies of $M\times\R$ identified in pairs along their
edges. A fundamental domain for the action is the regular polygonal prism $M\times I$ ($I$ a compact interval) illustrated in Figure \ref{Fig2}. We recover $X_m$ from this fundamental domain as follows. First identify
the top and bottom 
$m$-gons by a 2$(2\pi/m)$-twist to obtain a solid torus $N$.
(This corresponds to 
the quotient by $\phi^2$ above.) 
The vertical edges (of the form $\{v\}\times I$ where $v$ is a vertex of the 
$m$-gon) combine to give curves on the boundary of $N$ : a single circle $\xi_0$ in the case $m$ odd, and a pair of parallel circles $\xi_1,\xi_2$ if $m$ is even. The identification of these curves by $\phi|_E$ then either glues $\xi_0$ to itself by the antipodal map ($m$ odd) or glues $\xi_1$ to $\xi_2$ ($m$ even). The image of these curves shall be written $\xi$ in each case. The result of these identifications is $X_m$.

Note that $X_m$ inherits from $M\times{\R}$ a natural piecewise 
Euclidean metric, but that this metric depends on the two factors:
the area of $M$ and the translation length of $\phi$.

We choose a basepoint $v$ in $X_m$ and indicate, in Figure~\ref{Fig2}, a pair of closed geodesics in $X_m$, labeled $a$ and $b$, which pass through $v$. (These appear as geodesic paths in the universal cover, or in the prism $M\times I$.)    
The precise metric on $X_m$ shall now be chosen so that $a$ and $b$ both have length $1$, and then depends only on the acute angle $\delta$ between each of these closed geodesics and the vertical lines (see Figure~\ref{Fig2}).
We write $X_m(\delta)$ when we wish to make the metric explicit. Note that the metric may be chosen so that $\delta$ takes any value in the range $0<\delta<\frac{\pi}{2}$.

The closed geodesics $a$ and $b$ will be called {\it generator curves}.
They represent generators $a$ and $b$ for the fundamental group so that $\pi_1(X_m,v) = A(m;a,b)$. 
This may be seen by applying Van Kampen's Theorem. The element $x=ab$ is represented by a curve running once 
around the core of  the solid torus $N$, and $\Delta=(a,b)_m$ is represented by the curve $\xi$. Observe that
$$
[\xi_0] =x^m= \Delta^2 =z\ \hbox{ if $m$ odd, and } 
$$
$$
[\xi_1] =x^{m/2}= \Delta =z\ \hbox{ if $m$ even, }
$$

\noindent where $z=z_{a,b}$ generates the centre of $A(m)$.
     
\begin{figure}[ht]
\begin{center}
\includegraphics[width=3.2in]{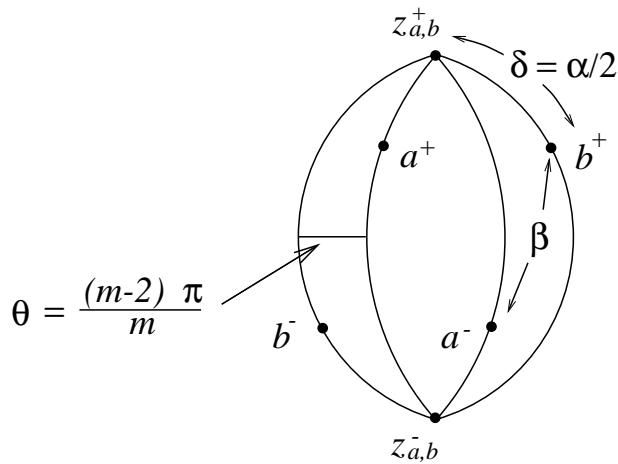}
\end{center}
\caption{The link of the vertex $v$ in the 3-complex $X_m$. }\label{Fig3}
\end{figure}

In Figure~\ref{Fig3} we see the link of $v$ in $X_m$. It is a piecewise spherical 
2-complex which can be described as the orthogonal join 
of a $0$-sphere with the union of two geodesic segments of length $\theta = \frac{(m-2)\pi}{m}$. We have indicated the points $a^{\pm}$ and $b^{\pm}$ 
corresponding to the geodesic loops based at $v$ 
which represent the generators $a$ and $b$ respectively.  
The $0$-sphere points correspond to the geodesic loop $\xi$ 
(based at $v$) which runs in the direction of the
the central element $z_{a,b}$.

It is easy to see, by checking the link condition, that $X_m$ is locally CAT(0). Its universal 
cover is a CAT(0) space which decomposes as a metric product of ${\R}$ 
with a 2-complex which is composed of infinitely many regular 
euclidean $m$-gons attached together along their vertices in a tree-like 
pattern. Each $m$-gon is attached to $m$ others, one along each of its 
$m$ vertices. Note that this 2-complex can be $A(m)$-equivariantly 
retracted down onto a Bass-Serre tree for $A(m)/Z=\Z_m\ast\Z_2$. 
In fact one can see in this way that $X_m$ has, as a $2$-spine, an equivariantly embedded copy of the CAT(0) 2-complex described by T. Brady and J. McCammond in \cite{BMcC}.

\subsection{Combining the building blocks}

We shall use the complexes  $X_m(\delta)$, $0<\delta<\frac{\pi}{2}$, as building blocks for constructing complexes 
for Artin groups with three or more generators. We adopt the convention that 
$X_2(\delta)$ is a 2-torus where the angle between geodesic representatives of the two fundamental group generators is $2\delta$.

Note that in the link of the vertex 
$v$ in the complex $X_m(\delta)$, as shown in Figure \ref{Fig3},
the points $a^\pm$ and $b^\pm$  are all a distance $\delta$
from the zero sphere $\{z_{a,b}^\pm\}$.

Define angles $\alpha$ and $\beta$ by  
$$
\alpha \; = \; d(a^+, b^+) \; = \; d(a^-,b^-) \; = \; 
2\delta
$$
and 
$$
\beta \; = \; d(a^+,b^-) \; =\; d(a^-,b^+)\,.
$$

When we come to combining the $X_m(\delta)$ spaces by gluing them along the $a$ and $b$ generator curves, the 
resulting complexes will have vertex links obtained 
by combining several of the links in Figure~\ref{Fig3} 
together along the $\{a^\pm\}$ and $\{b^\pm\}$ spheres. The key geometric 
observation in this section is that for large values of $m$ one can 
arrange that both $\alpha$ and 
$\beta$ are large enough that
the new links one obtains by such combinations will still be CAT(1). 

First we need a lemma describing the precise relationship between 
 $\alpha$, $\beta$ and $m$.

\begin{lemma}\label{trigeqn}
Let $m \geq 3$ be  an integer, and let 
$\theta = \frac{(m-2)\pi}{m}$. 
Then $\alpha$, $\beta$ and $\theta$ are related by the equation
$$
2\cos\beta \, + \, \left( 1+ \cos\theta \right)\cos\alpha 
\, + \, \left( 1-\cos\theta\right) \; = \; 0 \,.
$$
\end{lemma}

\Proof We have $d(a^-,z_{a,b}^-) = \delta=\alpha/2 = d(b^+,z_{a,b}^+)$, and 
we wish to compute $\beta =d(a^-,b^+)$. 
Realize the half of the link in Figure \ref{Fig3} which contains $a^-$ and $b^+$ 
as a subset of the unit sphere in ${\R}^3$, with $z_{a,b}^+$, 
$z_{a,b}^-$, $a^-$ and $b^+$ mapping 
to the unit vectors  $\langle 0,0,1\rangle$,  $\langle 0,0,-1\rangle$, 
$\langle \sin\delta, 0, -\cos\delta\rangle$ and 
$\langle \sin\delta\cos\theta, \sin\delta\sin\theta, 
\cos\delta \rangle$ respectively. 

Then 
\[
\begin{aligned}
\cos\beta\; &=\; 
\langle \sin\delta, 0, -\cos\delta\rangle \cdot 
\langle \sin\delta\cos\theta, \sin\delta\sin\theta, 
\cos\delta\rangle \\
 &= \; \sin^2\delta\cos\theta \, - \, \cos^2\delta\\
 &=\; \left( \frac{1-\cos\alpha}{2}\right)\cos\theta \, - \, 
\left(\frac{1+\cos\alpha}{2} \right)
\end{aligned}
\]
and the result follows on rearranging terms. 

\endproof

For example, in the case $\alpha = \beta$ we can 
solve the equation in Lemma~\ref{trigeqn} above to get  
$$
\cos \alpha \; =\; \frac{\cos\theta - 1}{\cos\theta + 3}\, .
$$
Table~1 shows how the $\alpha = \beta$ values depend on $m$. 
Angle measurements are in degrees. It will be 
useful in proving Theorem~\ref{3complexes} below. 

\begin{table}[p]
\begin{center}
\begin{tabular}{|p{0.1in}p{0.6in}p{1in}p{0.6in}p{1in}p{0.8in}|} 
\hline
 & & & & & \\
 & $m$ & $\theta=\frac{(m-2)\pi}{m}$ &$\cos\theta$ &  $\cos\alpha =\cos\beta$ 
&  $\alpha = \beta$\\
 & & & & & \\
\hline
& & & & & \\
& 3 & 60 & 0.5 & -1/7 & \ 98.213 \\
& & & & & \\
& 4 & 90 & 0 & -1/3 & 109.471 \\
& & & & & \\
& 5 & 108 & -0.309 & -0.486 & 119.107 \\
& & & & & \\
& 6 & 120 & -0.5 & -3/5 & 126.870 \\
& & & & & \\
& 7 & 128.571 & -0.623 & -0.683 & 133.090 \\
& & & & & \\
& 8 & 135 & -$\sqrt{2}/2$ & -0.745 & 138.118 \\
& & & & & \\
& 9 & 140 & -0.766 & -0.791 & 142.237 \\
& & & & & \\
& 10 & 144 & -0.809 & -0.826 & 145.656 \\
& & & & & \\
& 11 & 147.273 & -0.841 & -0.853 & 148.531 \\
& & & & & \\
& 12 & 150 & -0.866 & -0.874 & 150.978 \\
& & & & & \\
& 13 & 152.307 & -0.885 & -0.892 & 153.083 \\
& & & & & \\
& 18 &  160 & -0.940 & -0.941 & 160.298 \\
& & & & & \\
& 19 & 161.053 & -0.946 & -0.947 & 161.306 \\
& & & & & \\
& 21 & 162.857 & -0.9556 & -0.9565 & 163.046 \\
& & & & & \\
& 22 & 163.64 & -0.9595 & -0.9603 & 163.801 \\
& & & & & \\
& 43 & 171.628 & -0.9893 & -0.9894 & 171.650 \\
& & & & & \\
& 44 & 171.818 & -0.9898 & -0.9899  & 171.839  \\
& & & & & \\
\hline
\end{tabular}
\end{center}
\caption{Dependence of $\alpha = \beta$ values on $m$.}\label{Table1} 
\end{table}

\bigskip

\begin{thm}\label{3complexes}
All but finitely many of the 3-generator 
Artin groups of the form $A(m,n,2)$ with $3\leq m,n<\infty$
are the fundamental groups of compact locally CAT(0) 3-complexes.
\end{thm}

\Proof We construct a space by gluing together the building blocks for the $A(m)$, $A(n)$ and $A(2)$ as follows. Let
\[
X_{m,n,2}= X_m(\delta_1)\cup X_n(\delta_2) \cup X_2(\pi/2) / \sim
\]

where we identify all three basepoints to a point $v$ and identify the generator curves from each building block in pairs in the obvious fashion so that the fundamental group of $X_{m,n,2}$ is $A(m,n,2)$. 
We now compute the values of $m$ and $n$ for which $X_{m,n,2}$  is locally CAT(0). 

To apply the link condition we first observe that all
the 2-dimensional piecewise spherical links in $X$ are locally CAT(1), and so will be CAT(1) precisely when they 
have no closed geodesic loops of length less than $2\pi$. 

\begin{figure}[ht]
\begin{center}
\includegraphics[width=3.5in]{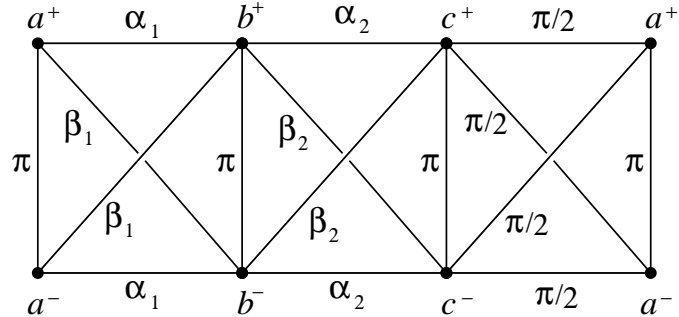}
\end{center}
\caption{Schematic picture of link of a vertex in the complex for $A(m,n,2)$.}\label{Fig4} 
\end{figure}

Suppose for concreteness that $A(m,n,2)$ has Artin generators $a$, $b$ and 
$c$ where $a$ and $b$ satisfy the relation of length $2m$, $b$ and $c$ 
satisfy the relation of length $2n$, and $a$ and $c$ commute. The link 
of the vertex in the complex $X_2(\pi/2)$ for $\langle a,c\rangle$ is a circle with 
two linking zero spheres labeled by $\{a^+, a^-\}$ and by $\{c^+, c^-\}$. 
The link of the vertex in the 3-complex for $\langle a,b \rangle$ (and 
likewise for $\langle b,c\rangle$) is as shown in Figure~\ref{Fig3}. These links 
combine together by identifying points with the same 
label from the set 
$\{a^+, a^-,  b^+, b^-,c^+, c^-\}$. 
The building block links are all CAT(1), so the only place to 
look for short loops is when one combines three links as above. 
We represent the combined link schematically as shown in Figure~\ref{Fig4}, where we use subscripts to distinguish the (possibly different) values for $\alpha$ and $\beta$ in each of $X_m(\delta_1)$ and $X_n(\delta_2)$.
There are six distances to measure between the points of any two of the three $0$-spheres $\{a^+, a^-\}$, $\{b^+, b^-\}$ and $\{c^+, c^-\}$. We record these distances schematically using a complete graph on 4 vertices.
The link shown in Figure \ref{Fig4} is the union of three such graphs.  

Working with  model 3-complexes  where $\alpha_i = \beta_i$ 
it is clear that one just needs to ensure that 
\[
\alpha_1 + \alpha_2 + \pi/2 \geq 2\pi 
\]
in order to get combined links which have no loops of length less than 
$2\pi$.  There are a few cases to consider. 

\medskip
\noi
{\it Case $A(m,3,2)$}. In this case $\alpha_2 = \beta_2 = 98.21$ from Table~1. 
We see that $m$ must be chosen large enough so that $\alpha_1 = \beta_1$ 
is greater than $360 - (98.21 + 90) = 171.79$. 
Table~1 shows that this is achieved for $m \geq 44$. 

\medskip
\noi
{\it Case $A(m,4,2)$}. In this case $\alpha_2 = \beta_2 = 109.5$ from Table~1. 
We see that $m$ must be chosen large enough so that $\alpha_1 = \beta_1$
is greater than $360 - (109.5 + 90) = 160.5$. 
Table~1 shows that this is achieved for $m \geq 19$. 

\medskip
\noi
{\it Case $A(m,5,2)$}. In this case $\alpha_2 = \beta_2 = 119.1$ from Table~1. 
We see that $m$ must be chosen large enough so that $\alpha_1 = \beta_1$
is greater than $360 - (119.1 + 90) = 150.9$. 
Table~1 shows that this is achieved for $m \geq 12$.

\medskip
\noi
{\it Case $A(m,6,2)$}. In this case $\alpha_2 = \beta_2 = 126.8$ from Table~1. 
We see that $m$ must be chosen large enough so that $\alpha_1 = \beta_1$
is greater than $360 - (126.9 + 90) = 143.1$. 
Table~1 shows that this is achieved for $m \geq 10$.

\medskip
\noi
{\it Case $A(m,7,2)$}. In this case $\alpha_2 = \beta_2 = 133.1$ from Table~1. 
We see that $m$ must be chosen large enough so that $\alpha_1 = \beta_1$
is greater than $360 - (133.1 + 90) = 136.9$. 
Table~1 shows that this is achieved for $m \geq 8$.

Clearly all groups of the form $A(m,n,2)$ with $m,n \geq 8$ will give 
3-complexes with CAT(1) links. Hence the theorem follows. 
 \endproof

We can combine this result with previously known results 
to get a more complete 
picture of all 3-generator Artin groups.  

\begin{corollary}\label{3gens}
All but finitely many of the 3-generator Artin groups 
are the fundamental groups of compact locally CAT(0) complexes of dimension 
at most three. 
\end{corollary}

\Proof The Artin groups of the form $A(m,n,p)$ where $m,n,p \geq 3$,
and the three generator Artin groups with fewer than three relations 
all admit 2-dimensional CAT(0) structures by the work of \cite{BMcC}. The 
Artin groups of the form $A(2,2,p)$ $p\geq 2$ all admit 3-dimensional 
CAT(0) structures by \cite{Bra}. Theorem~\ref{3complexes} above covers the 
remaining cases ($A(m,n,2)$ where $m,n \geq 3$). 
\endproof

It is worth noting that we could combine the $A(m)$ building blocks 
to produce compact, locally CAT(0) 3-dimensional Eilenberg Mac Lane 
complexes for all but finitely many of the 3-generator Artin groups of 
\emph{large type}; where $m$, $n$ and $p$ are all at least 3. This is made explicit 
in Theorem \ref{manygens}(1) below. However, in the 3-generator case above, it is more 
efficient to appeal to the neat 2-dimensional construction of \cite{BMcC}.   

Note that Theorem~\ref{3complexes} and Corollary~\ref{3gens} exclude 
precisely 65 infinite type Artin groups $A(m,n,2)$ in addition to the  
finite types $A(3,3,2)$, $A(4,3,2)$ and $A(5,3,2)$ not considered in the proofs. We note that 3-dimensional 
locally CAT(0) Eilenberg-Mac Lane complexes have been found for the finite 
type three generator Artin groups \cite{Bra}. The CAT(0) dimension 
is unknown for the remaining 65 examples.

Our 3-dimensional 
building blocks can also be combined with some limited success in the 
case of Artin groups with more than 3 generators. Recall that an Artin group 
is specified by a finite set $\{ a_1, \ldots , a_n\}$ of generators, and 
each pair $\{a_i, a_j\}$ of generators satisfies a relation of the form 
\[
(a_i,a_j)_{m_{ij}} \; = \; (a_j,a_i)_{m_{ij}}
\]
where the $m_{ij}$ belong to $\{2,3,4, \ldots \} \cup \{\infty\}$. 
In the case $m_{ij} = \infty$ we have no explicit 
relation between $a_i$ and $a_j$. 
We call the $m_{ij}$ {\it relator indices}. 

\begin{thm}\label{manygens}
The following types of Artin groups all admit compact, 
3-dimensional, locally CAT(0) Eilenberg Mac Lane complexes. 
\begin{enumerate}

\item Artin groups with the property that no 
triple of relator indices of the form 
$\{m_{ij}, m_{jk},m_{ki} \}$ appears in the following list
\begin{tabbing}
\hskip1cm\=\hskip5cm\=\hskip5cm\=\kill
\>$\{m,2,2\}$, $2\leq m<\infty$; \>$\{m,3,3\}$, $3\leq m <22$; \>$\{m,4,4\}$, $4\leq m <9$;\\
\>$\{m,3,2\}$, $3\leq m <44$;    \>$\{m,4,3\}$, $4\leq m <13$; \>$\{5,5,4\}$, $\{6,5,4\}$; \\
\>$\{m,4,2\}$, $4\leq m <19$;    \>$\{m,5,3\}$, $5\leq m <10$; \>\\
\>$\{m,5,2\}$, $5\leq m <12$;    \>$\{6,6,3\}$, $\{7,6,3\}$; \>$\{5,5,5\}$. \\
\>$\{m,6,2\}$, $6\leq m <10$; \>\> \\
\>$\{7,7,2\}$\>\>\\
\end{tabbing}

\item Artin groups with the property that all the relator indices $m_{ij}$ 
are  at least $4$, and such that if $\{m_{ij}, m_{jk},m_{ki}\}$ is a triple
of finite indices, then $\{m_{ij}, m_{jk},m_{ki}\}\neq\{4,4,4\}$ 
\end{enumerate}
\end{thm}

\Remark Part 1 of this Theorem includes the class of \emph{triangle free} Artin groups, those with the property that, for each triple of relator indices of the form $\{m_{ij}, m_{jk},m_{ki} \}$, at least one of the indices is infinite. Triangle free Artin groups were shown to have 2-dimensional locally CAT(0) Eilenberg-MacLane spaces in \cite{BMcC}.
However, we note that, for the most part, the Artin groups listed in
Theorem \ref{manygens} are not covered by the 2-dimensional results of T.~Brady and J.~McCammond.

\Proof The proof proceeds in the same way as the proof of Theorem~\ref{3complexes}. 
We glue the model complexes for the dihedral type Artin groups together, and 
then check to see if the links in the resulting complex are CAT(1). As 
before, this reduces to checking that there are no closed loops of length 
less than $2\pi$ in the combined links. Again this reduces to checking that 
there are no short loops in the corresponding schematic picture of the link. 

Since all $\alpha$ and $\beta$ are at least $\pi/2$,  we  only have to 
consider triples of finite indices of the form  
$\{m_{ij}, m_{jk},m_{ki}\}$ when looking for short loops in the 
schematic graph.

In the case of Artin groups of type~(1) we use model spaces where 
$\alpha = \beta$. By employing Table \ref{Table1} and the approach used in the proof of Theorem \ref{3complexes} we find that there is a short loop in the schematic graph associated to a triple $\{m_{ij}, m_{jk},m_{ki}\}$ of
finite indices if and only if it is in the given list.

In the case of Artin groups of type~(2) we generally choose $\alpha\neq\beta$. We choose the model space $X_4(\delta)$ for $A(4)$ by letting $\al=\al_4 = 163^\circ$. Then using Lemma \ref{trigeqn} (with $m=4$ we have $\cos\theta =0$) we have $\be=\be_4 = 91.25^\circ > 91^\circ$. 

We remark that, regardless of the choice of model space $X_m(\delta)$, $0< \delta < 90^\circ$, we will always have $\be > \theta$.
(Note that $0<\theta <180^\circ$ for all $m$. Fix $\theta$ in this range. From the equation 
$\cos\beta = \sin^2\delta\cos\theta -\cos^2\delta$ in the proof of Lemma \ref{trigeqn}, we have 
that $\cos\be =\cos\theta$ if and only if $\delta =\pm 90^\circ$, while at $\delta =0^\circ$, 
$\cos\be = -1$. So, by continuity, $\cos\be < \cos\theta$ for all $\delta$ in the range 
$0< \delta < 90^\circ$. Hence $\be >\theta$.) 
Thus, when $m\geq 5$ we will always have $\beta > 108^\circ$, leaving us free to choose $\alpha$ as close to $180^\circ$ as we like.
We choose the spaces $X_m(\delta)$ for $m\geq 5$ by putting $\alpha=\al_{\geq 5}=179^\circ$ ($\delta=89.5^\circ$), so ensuring that
$\al_{\geq 5}+2\be_4>360^\circ$.

Again we consider the triples $\{m_{ij}, m_{jk},m_{ki}\}$ of finite indices, 
all at least 4. Clearly the sum of three $\alpha$ edges will always be 
greater than $3(163^\circ)>360^\circ$. Otherwise, we must consider the sum 
of one $\alpha$ edge with two $\beta$ edges (one edge coming from the 
schematic link component associated with each of the three relator indices 
in $\{m_{ij}, m_{jk},m_{ki}\}$). If one of the $\beta$ edges comes from an index greater than $5$ (i.e: its length is greater than $108^\circ$), then the sum 
of the three will be greater than
$163^\circ + 91^\circ + 108^\circ = 362^\circ$.
This leaves the case where both $\beta$ edges are short (coming from relator 
indices equal to $4$).  We have already seen that $\al_{\geq 
5}+2\be_4>360^\circ$, leaving only the possibility that $\{m_{ij}, 
m_{jk},m_{ki}\} =\{4,4,4\}$ which is excluded by the hypothesis. We note, 
unfortunately, that there is no possible choice of $\al_4$ which allows 
$\al_4+2\be_4$ to exceed $360^\circ$. 
\endproof

\section{Related results and further questions}\label{Sect5}




In \cite{Bri2} Bridson presents a group which has geometric dimension $2$ and CAT(0) dimension 3 but which, remarkably, contains an index 2 subgroup with CAT(0) dimension 2. This prompts the following question: 
 
\begin{question}
Do the groups $A(m,n,2)$, with $m,n$ odd, contain finite index subgroups with CAT(0) dimension $2$? 
\end{question}

This seems plausible in view of P. Hanham's results in \cite{Han}, 
which indicate that the proof of Theorem~A relies critically on the particular 
structure of minsets  that is forced by the presence 
of  ``odd'' Artin relators.  This is very specific geometric information about 
specific elements of the group -- information that can be lost after
passing to a finite index subgroup.

We remarked  in Section~\ref{Sect3} after the proof of  
Theorem~\ref{3complexes} and Corollary~\ref{3gens} that our 
3-dimensional constructions did not give CAT(0) structures for  
65 Artin groups of type $A(m,n,2)$ which are not finite type.
Following the results of \cite{Han} cited above, this list can be reduced to
the following $24$ groups:  
$A(3,2n+1,2)$ for $3\leq n \leq 21$, 
$A(2m+1,5,2)$ for $2\leq m \leq 5$,
$A(7,7,2)$.

\begin{question} What is the CAT(0) dimension of each of the 24 
Artin groups listed above? 
\end{question}

It is easily enough seen that taking a direct product of a group with
$\Z^n$ increases both the geometric and CAT(0) dimensions by precisely
$n$. For the CAT(0) dimension one uses the Flat Torus Theorem (Theorem
\ref{FTT}) and a result of Morita \cite{Mor} which states that the
covering dimension of $X\times Y$ is the sum of the dimension of $X$
and $Y$ in the case where $X$ is a locally finite polyhedron and $Y$ any
metric space. This observation yields examples of groups of arbitrary
finite geometric dimension which exhibit a gap of 1 between the geometric
and CAT(0) dimensions, by simply taking products of the existing
2-dimensional examples with $\Z^n$. In principal, as was pointed
out in \cite{Bri2}, one can also obtain
arbitrarily large (finite) gaps by taking products of existing examples
with themselves. However in practice one may need to introduce extra
hypotheses to achieve this. For example, if one defines ``CAT(0)
dimension'' over the class of complete CAT(0) spaces with the geodesic
extension property and using Kleiner's dimension for CAT(0) spaces
\cite{Kl}, then this dimension may be shown to be additive with
respect to products of centreless groups. In contrast to the
topological covering dimension, Kleiner's dimension is additive
with respect to products; the remaining hypotheses allow one to use the
splitting theorem II.6.21 of \cite{BH}.

One motivation for the $2$-dimensional structures found in \cite{BMcC} was to show that the associated Artin groups are biautomatic. 

\begin{question}
Can one use the 3-dimensional CAT(0) structures to show that the Artin groups listed in Theorem \ref{manygens} are biautomatic? 
\end{question}

The groups of Theorem \ref{manygens}, part (2), are contained in the class
of extra-large type Artin groups (all $m_{i,j}\geq 4$) which have already been shown to be biautomatic by D. Peifer \cite{P}. According to Corollary III.$\Gamma$.4.8 of \cite{BH}, any group which acts properly and cocompactly on a CAT(0) space is semi-hyperbolic (asynchronously bi-combable) and hence has solvable word and conjugacy problems.

There are other variations of the notion of CAT(0) dimension and many 
interesting questions to consider. There is, for example, the question mentioned in the Introduction
as to what happens if one removes the semi-simplicity assumption.   
One can also define the CAT(-1) dimension of a group exactly as we 
defined CAT(0) dimension in this paper by replacing the occurrences of 
CAT(0) with CAT(-1). 

\begin{question}\label{hyp1}
Are there torsion free Gromov hyperbolic groups with distinct CAT(0) and CAT(-1) dimensions? 
\end{question}

The answer is yes if one restricts attention to actions on \emph{proper} CAT(-1) spaces. A first example is given in \cite{BC2}. This example arises as a finite index subgroup of an example with torsion due to M. Kapovich \cite{Ka} (see also \cite{BC}). 

It is unknown if Gromov hyperbolic groups are all CAT(0) or CAT(-1), that is whether
they act properly and cocompactly on a finite dimensional CAT(0), respectively CAT(-1),
complex. Perhaps the following question is more manageable. 
Note that we do not require cocompactness in the definitions 
of CAT(0) and  CAT(-1) dimensions. 

\begin{question}\label{hyp2}
Do all torsion free Gromov hyperbolic groups have 
finite CAT(0) or CAT(-1) dimensions? For a positive result
one would even be happy to drop the semi-simplicity condition.
\end{question}

We know that the Rips complex gives a finite Eilenberg-MacLane space 
for torsion free Gromov hyperbolic groups, so they do have finite geometric 
dimension.  

Of course our definition of CAT(0) dimension does not restrict us to looking at torsion free groups. Question \ref{hyp2}, for example, is still of considerable interest in the case of hyperbolic groups with torsion. Some interesting examples of hyperbolic groups with torsion are considered in \cite{BC}. Here the CAT(0) dimension is compared with the minimal dimension of a proper $G$-CW-complex, an approach which is also in the spirit of finding the  minimum dimension of a nonpositively curved $\underbar{E}(G)$.

\end{document}